\documentclass[10pt,notitlepage,reqno]{amsart}

\usepackage{amsmath}
\usepackage{amsthm}

\DeclareMathOperator{\ev}{ev}
\DeclareMathOperator{\Ind}{Ind}

\usepackage{amsfonts}
\usepackage{scalerel}
\usepackage{amssymb,amscd,epsf,verbatim,mathtools,framed}
\usepackage{mathrsfs}
\usepackage{graphicx}
\usepackage{placeins}
\usepackage{needspace}
\usepackage{latexsym}
\usepackage{lscape}
\usepackage[colorlinks=true]{hyperref}
\hypersetup{colorlinks, citecolor=blue, filecolor=black, linkcolor=purple, urlcolor=violet}
\usepackage{epstopdf}
\usepackage{pdflscape}
\usepackage{enumitem}

\usepackage{tikz}
\usepackage[all]{xy}
\usetikzlibrary{calc}
\usetikzlibrary{matrix,arrows,decorations.pathmorphing}
\usepackage{tabularx}
\usepackage{tikz-cd}
\usepackage{color}
\usepackage{multirow}
\usepackage{scalefnt}
\usepackage{fancyhdr}
\usepackage[margin=1.2in]{geometry}
\usepackage[T1]{fontenc}

\theoremstyle{definition}

\newtheorem{thm}{Theorem}[section]
\newtheorem{example}[thm]{Example}
\newtheorem{defn}[thm]{Definition}
\newtheorem{prop}[thm]{Proposition}

\newtheorem{lem}[thm]{Lemma}
\newtheorem{cor}[thm]{Corollary}
\newtheorem{construction}[thm]{Construction}

\newtheorem{rmk}[thm]{Remark}

\makeatother
\numberwithin{equation}{section}

\newcommand{\cA}{{\mathcal A}}
\newcommand{\cB}{{\mathcal B}}
\newcommand{\cC}{{\mathcal C}}
\newcommand{\cD}{{\mathcal D}}
\newcommand{\cE}{{\mathcal E}}
\newcommand{\cF}{{\mathcal F}}
\newcommand{\cG}{{\mathcal G}}
\newcommand{\cH}{{\mathcal H}}
\newcommand{\cI}{{\mathcal I}}
\newcommand{\cJ}{{\mathcal J}}

\newcommand{\cL}{{\mathcal L}}
\newcommand{\cM}{{\mathcal M}}

\newcommand{\cO}{{\mathcal O}}
\newcommand{\cP}{{\mathcal P}}

\newcommand{\cS}{{\mathcal S}}
\newcommand{\cT}{{\mathcal T}}

\newcommand{\cV}{{\mathcal V}}
\newcommand{\cW}{{\mathcal W}}
\newcommand{\cX}{{\mathcal X}}

\newcommand{\cZ}{{\mathcal Z}}

\newcommand{\kS}{{\mathfrak S}}

\newcommand{\kg}{{\mathfrak g}}

\newcommand{\kl}{{\mathfrak l}}

\newcommand{\kh}{{\mathfrak h}}

\newcommand{\bA}{{\mathbb A}}
\newcommand{\bB}{{\mathbb B}}
\newcommand{\bC}{{\mathbb C}}

\newcommand{\bE}{{\mathbb E}}
\newcommand{\bF}{{\mathbb F}}
\newcommand{\bG}{{\mathbb G}}

\newcommand{\bL}{{\mathbb L}}
\newcommand{\bM}{{\mathbb M}}
\newcommand{\bN}{{\mathbb N}}
\newcommand{\bP}{{\mathbb P}}
\newcommand{\bQ}{{\mathbb Q}}

\newcommand{\bT}{{\mathbb T}}

\newcommand{\bV}{{\mathbb V}}
\newcommand{\bW}{{\mathbb W}}

\newcommand{\bZ}{{\mathbb Z}}

\newcommand{\mH}{{\mathrm H}}
\newcommand{\mR}{{\mathrm R}}

\newcommand{\et}{\mathrm{\acute{e}t}}   

\newcommand{\dR}{\mathrm{dR}}

\newcommand{\st}{{\mathrm{st}}}

\newcommand{\ad}{\mathrm{ad}}

\newcommand{\an}{\mathrm{an}}

\newcommand{\ab}{\mathrm{ab}}

\newcommand{\id}{\mathrm{id}}

\newcommand{\pr}{\mathrm{pr}}
\DeclareMathOperator{\Gal}{Gal}

\newcommand{\HH}{\mathrm{HH}}
\newcommand{\nil}{\mathrm{nil}}

\DeclareMathOperator{\Spec}{Spec}    
\DeclareMathOperator{\Spa}{Spa}
\DeclareMathOperator{\Ker}{Ker}

\DeclareMathOperator{\tr}{tr}

\DeclareMathOperator{\End}{End}

\DeclareMathOperator{\Loc}{Loc}
\DeclareMathOperator{\Hom}{Hom}
\DeclareMathOperator{\Ext}{Ext}
\DeclareMathOperator{\sEnd}{\underline{End}}
\DeclareMathOperator{\sHom}{\underline{Hom}}
\DeclareMathOperator{\sExt}{\underline{Ext}}
\DeclareMathOperator{\sSpec}{\underline{Spec}}
\DeclareMathOperator{\MIC}{MIC}

\DeclareMathOperator{\Lie}{Lie}
\DeclareMathOperator{\Rep}{Rep}

\DeclareMathOperator{\rk}{rank}
\DeclareMathOperator{\GL}{GL}

\DeclareMathOperator{\Uni}{Uni}
\DeclareMathOperator{\uni}{uni}

\renewcommand{\leq}{\leqslant}
\renewcommand{\geq}{\geqslant}

\newcommand{\ve}{\varepsilon}    

\newcommand{\vp}{\varphi}
\newcommand{\ra}{\rightarrow}

\newcommand{\hra}{\hookrightarrow}
\newcommand{\xra}{\xrightarrow}

\newcommand{\ul}{\underline}
\newcommand{\barK}{{\overline{K}}}

\newcommand{\OK}{{\cO_K}}
\newcommand{\Qp}{{\bQ_p}}
\newcommand{\Zp}{{\bZ_p}}

\title{Rigidity of good reduction in families of curves}
\author{Zhefan Duan}
\date{}

\begin{document}

\begin{abstract}
Consider a family $X\rightarrow S$ of geometrically connected smooth projective curves where $S$ is a geometrically connected smooth rigid analytic variety with good reduction. Assume that every fiber has genus $\geq 1$ and has semistable reduction. We prove the following rigidity result: if the fiber $X_{s_0}$ over some classical point $s_0$ has good reduction, then the family has good reduction everywhere. This is a non-abelian version of the rigidity theorem of Diao--Yao. We also prove an analogue of the rigidity result when the base $S$ has semistable reduction. The proofs reply on the construction of a canonical $\Zp$-local system on $S$ that detects the good reduction property of the fibers.
\end{abstract}

\maketitle
\thispagestyle{empty}
\tableofcontents

\section{Introduction}

\subsection{Main results on rigidity of good reduction}
\noindent
\vspace{0.05in}

\noindent
Throughout the article, let $p$ be a prime number and let $K$ be a $p$-adic field, that is, a complete discretely-valued field of mixed characteristic $(0,p)$ with a perfect residue field. Let $\mathcal{O}_K$ be the ring of integers and $k$ be the residue field. Let $\Gal_K= \Gal(\overline K/K)$ denote the absolute Galois group of $K$.

We study smooth proper families of curves $q \colon X\ra S$ where $S$ is a geometrically connected smooth rigid analytic variety over $K$. For any classical point $s$ of $S$, we use $\kappa(s)$ to denote its residue field and use $X_s$ to denote the fiber over $s$. The main result is the following.

\begin{thm}\label{thm: introduction rigidity}
Let $\kS$ be a smooth $p$-adic formal scheme over $\mathcal{O}_K$ with geometrically connected special fiber. Let $S=\kS_{\eta}$ be the rigid analytic generic fiber and let $q \colon X\ra S$ be a smooth proper family of geometrically connected curves of genus $\geq 1$. Suppose that $X_s$ has semistable reduction over $\kappa(s)$ at every classical point $s\in S$. If there is a classical point $s_0$ of $S$ such that $X_{s_0}$ has good reduction over $\kappa(s_0)$, then $X_s$ has good reduction over $\kappa(s)$ for every classical point $s$ of $S$.
\end{thm}

Theorem~\ref{thm: introduction rigidity} indicates that the property of admitting good reduction is ``rigid'' in an analytic family of smooth projective curves, provided the base $S$ has good reduction and the fibers already have semistable reduction. This can be viewed as a non-abelian incarnation of the rigidity theorem of Diao--Yao \cite[Theorem 1.1]{DY}.

\begin{rmk}
    The good reduction hypothesis on the base in Theorem~\ref{thm: introduction rigidity} is crucial; see \S \ref{subsec: failure of rigidity} for a counterexample when the base is a closed annulus, which does not have good reduction. 
\end{rmk}

When the base has only semistable reduction, we also prove an analogue of this rigidity result.
\begin{thm}\label{thm: introduction rigidity(semistable reduction)}
    Let $\kS$ be a semistable $p$-adic formal scheme over $\mathcal{O}_K$ with geometrically connected special fiber. Let $S=\kS_{\eta}$ be the rigid analytic generic fiber and let $q \colon X\ra S$ be a smooth proper family of geometrically connected curves of genus $\geq 1$. Suppose that $X_s$ has semistable reduction over $\kappa(s)$ at every classical point $s\in S$. 
    If for every irreducible component $T$ of $\kS_k$,  there is a classical point $s_T$ of $\mathrm{sp}^{-1}(T\cap \kS_k^{\mathrm{sm}})$\footnote{
    Here, $\mathrm{sp}\colon S \ra \kS_k$ is the specialization map and $\kS_k^{\mathrm{sm}}$ denotes the smooth locus inside $\kS_k$.
    } such that $X_{s_T}$ has good reduction over $\kappa(s_T)$, then $X_s$ has good reduction over $\kappa(s)$ for every classical point $s$ of $S$.
\end{thm}

\begin{rmk}
When the genus is $1$, the good reduction property can be detected by looking at the first \'etale cohomology, and hence Theorem~\ref{thm: introduction rigidity} and Theorem~\ref{thm: introduction rigidity(semistable reduction)} follow immediately from the main theorems of \cite{DY}. Therefore, we may assume that the genus is at least $2$.
\end{rmk}

\begin{rmk} 
    In a different direction, one may consider the following algebraic situation.
    Assume that $p> 2g+1$. Let $\cS$ be a geometrically connected smooth proper scheme over $\mathcal{O}_K$ and let $S$ be its generic fiber. Let $q \colon X\ra S$ be a smooth proper family of geometrically connected curves. Suppose that for some closed point $s_0$ of $S$ with $\kappa(s_0)$ unramified over $K$, the fiber $X_{s_0}$ has good reduction. Using the $\ell$-adic method developed in \cite{CT}, one can prove that $X_s$ has good reduction for every closed point $s$ of $S$, without assuming that the fibers already have semistable reduction. More precisely, when $p> 2g+1$, we can choose a prime $\ell\neq p$ such that $p \nmid |\mathrm{GL}_{2g}(\bZ_\ell)|$.
    Then the outer pro-$\ell$ representation (cf. \cite[\S 6.1]{CT})
    \[
       \rho_\ell \colon \pi_1^{\et}(X) \ra \mathrm{Out}(\pi_1^\et(X_{\bar \eta})^{\{\ell\}})
    \]
    is tamely ramified. Therefore, the statement follows from \cite[Remark 24]{CT} and the proof of \cite[Corollary 4]{CT}, which uses Oda's criterion \cite{Oda1995}. This is closely related to Question 1.16 raised in \cite{DY}.
\end{rmk}

In the remainder of this introduction, we describe the key ingredients in the proofs of Theorem~\ref{thm: introduction rigidity} and Theorem~\ref{thm: introduction rigidity(semistable reduction)}, and then outline the overall proof strategy.

\vspace{0.05in}
\subsection{Rigidity of crystallinity}
\noindent
\vspace{0.05in}

\noindent
The first key ingredient is the rigidity theorems of Diao and Yao on crystallinity of $p$-adic local systems.

\begin{thm}[Diao--Yao, {\cite[Theorems 1.1]{DY}}]
\label{thm: introduction crystalline rigidity}
Let $\kS$ be a smooth $p$-adic formal scheme over $\mathcal{O}_K$ with geometrically connected special fiber. Let $S=\kS_{\eta}$ be the rigid analytic generic fiber and let $\bL$ be a semistable $\Zp$-local system on $S$. If $\bL|_{s_0}$ is crystalline at some classical point $s_0$ of $S$, then $\bL|_s$ is crystalline at every classical point $s$ of $S$.
\end{thm}

\begin{thm}[Diao--Yao, {\cite[Theorems 3.11]{DY}}]
\label{thm: introduction crystalline rigidity(semistable reduction)}
Let $\kS$ be a semistable $p$-adic formal scheme over $\mathcal{O}_K$ with geometrically connected special fiber. Let $S=\kS_{\eta}$ be the rigid analytic generic fiber and let $\bL$ be a semistable $\Zp$-local system on $S$. 
If for every irreducible component $T$ of $\kS_k$,  there is a classical point $s_T$ of $\mathrm{sp}^{-1}(T\cap \kS_k^{\mathrm{sm}})$ such that $\bL|_{s_T}$ is crystalline at $s_T$, then $\bL|_s$ is crystalline at every classical point $s$ of $S$.
\end{thm}

In fact, once we construct a suitable local system that detects the good reduction property of the fibers of $q\colon X\rightarrow S$, our main theorems follow immediately from Theorem~\ref{thm: introduction crystalline rigidity} and Theorem~\ref{thm: introduction crystalline rigidity(semistable reduction)}, respectively.

\vspace{0.05in}
\subsection{The non-abelian criterion of Andreatta--Iovita--Kim}
\noindent
\vspace{0.05in}

\noindent
A second key ingredient is the $p$-adic criterion for good reduction of curves due to Andreatta--Iovita--Kim in \cite{AIK}. As already pointed out in \emph{loc.\ cit.}, a smooth projective curve over $K$ may fail to have good reduction even if its Jacobian has good reduction. In particular, any criterion for good reduction of curves has to be non-abelian in nature.

To state the result of Andreatta--Iovita--Kim, we consider a geometrically connected smooth projective curve $Z$ over $K$ of genus $g\geq 2$. Assume that $Z$ has semistable reduction over $\mathcal{O}_K$. Let $\Uni_{\bQ_p}^{\leq n}(Z_{\barK})$ denote the category of unipotent $\Qp$-local systems on $Z_{\barK}$ of unipotency index at most $n$ (cf. \cite[\S 3.2]{AIK}). Moreover, pick a $K$-point $b\in Z(K)$ and a geometric point $\bar{b}$ above it. Consider the category $\Uni_{\bQ_p}^{\leq n}(Z_{\barK}, \bar{b})$ consisting of objects $(\bL,v)$ where $\bL \in \Uni_{\bQ_p}^{\leq n}(Z_{\barK})$ and $v\in \bL_{\bar{b}}$.

\begin{thm}[Andreatta--Iovita--Kim; {\cite[Proposition 3.5, Corollary 3.7, and Theorem 1.9]{AIK}}]
\label{thm: introduction AIK}
Let $Z$, $b$, and $\bar{b}$ be as above.
\begin{enumerate}
 \item For each $n\geq 1$, there exists a universal object $(\bF_n,u_n)$ in $\Uni_{\bQ_p}^{\leq n}(Z_{\barK}, \bar{b})$ in the sense of \cite[Definition 3.1]{AIK}. Moreover, $\bF_n$ uniquely descends to an \'etale $\Qp$-local system on $Z$, which we still denote by $\bF_n$ by a slight abuse of notation.
 \item The curve $Z$ has good reduction over $K$ if and only if the stalk $\bF_{4, \bar{b}}$ is crystalline as a $p$-adic Galois representation, equivalently, if and only if $\bF_{n, \bar{b}}$ is crystalline for every $n\geq1$.
\end{enumerate}
\end{thm}

As we shall see in a moment, the proof of our main theorems involves constructing, in the relative setting, a similar collection of universal unipotent local systems. Building on this, we are able to produce a local system that detects the good reduction property of the curves on fibers, at least in the special case where the family admits a section. This allows us to apply Theorem~\ref{thm: introduction crystalline rigidity} and Theorem~\ref{thm: introduction crystalline rigidity(semistable reduction)} and conclude.

\vspace{0.05in}
\subsection{Proof of the main theorems in the presence of a section}\label{subsection: intro_proof_of_main_thm_section}
\noindent
\vspace{0.05in}

\noindent
Let us first demonstrate the strategy for proving Theorem~\ref{thm: introduction rigidity} and Theorem~\ref{thm: introduction rigidity(semistable reduction)} in the special case where the family of curves (of genus $\geq 2$) $q: X\rightarrow S$ admits a section $b:S\rightarrow X$.

The first step is to generalize the construction in \cite{AIK} to the relative setting, that is, we construct a collection of universal local systems on $X$ that are ``unipotent relative to $S$''. More precisely, let $\Uni_{\Zp}^{\leq n}(q)$ denote the category of \'etale $\Zp$-local systems on $X$ that are \emph{relatively unipotent} with respect to $q$ with unipotency index $\leq n$ (see Definition~\ref{def: relative unipotence}). Also consider the category $\Uni_{\Zp}^{\leq n}(q,b)$ consisting of pairs $(\bL, v)$ with $\bL\in \Uni_{\Zp}^{\leq n}(q)$, together with a \emph{marking}
\[
 v \in \Hom(\ul\Zp,b^*\bL)= \mH^0_\et(S,b^*\bL).
\]
In fact, mimicking the construction in \cite{AIK}, we are able to produce a universal object $(\bE_n,u_n)$ in $\Uni_{\Zp}^{\leq n}(q,b)$ for each $n\geq 1$ (cf. Theorem~\ref{thm: existence of relative universal system}). In the case where $S$ is a point, $\bE_n[\frac 1p]$ coincides with the $\Qp$-local system $\bF_n$ in Theorem~\ref{thm: introduction AIK} (1). Moreover, for every $\bL\in\Uni_{\Zp}^{\leq n}(q)$, evaluation at $u_n$ induces an isomorphism
\begin{equation}\label{eq: intro universality}
 \mathrm{ev}_{u_n}\colon  q_*\sHom(\bE_n,\bL)\xra[]{\sim}b^*\bL. 
\end{equation}

We claim that the \'etale $\Zp$-local system $\bT\coloneqq b^*\bE_4$ on $S$ detects the good reduction property of the family of curves. To see this, pick any classical point $s$ on $S$ and consider the pullback diagram
\[
\begin{tikzcd}[column sep=1.6cm, row sep=1.2cm]
X_s \arrow[r,"\iota_s"] \arrow[d,"q_s"]
& X \arrow[d,"q"']\\
s \arrow[u, bend left=40pt, "b_s"{pos=0.5, yshift=6pt}]
  \arrow[r]
& S \arrow[u, bend right=40pt, "b"']
\end{tikzcd}
\]
where the section $b$ induces a section $b_s$ of $q_s$ (i.e., a $\kappa(s)$-point on $X_s$). Pick a geometric point $\bar{b}_s$ above $b_s$. For each $s$, let $\bF_n(s)$ denote the \'etale $\Qp$-local system on $X_s$ with respect to $b_s$ and $\bar{b}_s$ as in Theorem~\ref{thm: introduction AIK} (1). The base change property proved in Proposition~\ref{prop: base change universal system} yields an isomorphism
\[
 \iota_s^*\bE_n[\frac 1p]\cong \bF_n(s).
\]
In particular, we have $\bT_{{\bar s}}[\frac 1p]\cong \bF_{4}(s)_{\bar{b}_s}$. By Theorem~\ref{thm: introduction AIK}, if $X_s$ has semistable reduction over $\kappa(s)$, then $X_s$ has good reduction if and only if $\bT_{\bar{s}}$ is crystalline. Finally, Theorem~\ref{thm: introduction crystalline rigidity} and Theorem~\ref{thm: introduction crystalline rigidity(semistable reduction)} are precisely what we need to conclude the proofs of Theorem~\ref{thm: introduction rigidity} and Theorem~\ref{thm: introduction rigidity(semistable reduction)}, respectively.

\vspace{0.05in}
\subsection{Proof of the main theorem: general case}\label{subsec: intro_proof_main_thm_general}
\noindent
\vspace{0.05in}

\noindent
In the general setting, the family $q:X\rightarrow S$ might not have a section. The main idea is to consider the natural projection $\pr_1:X\times_SX\ra X$ onto the first component and consider the diagonal map $\Delta:X\ra X\times_SX$. Then $(\pr_1, \Delta)$ produces a smooth proper family of geometrically connected curves with a section, so that we can apply the construction in \S \ref{subsection: intro_proof_of_main_thm_section}. For each $n\geq 1$, write $\bE_{n,X}^\Delta$ for the universal unipotent $\Zp$-local system in $\Uni_{\Zp}^{\leq n}(\pr_1,\Delta)$, and set
\[
 \bA_{n,X}^\Delta\coloneqq \Delta^*\bE_{n,X}^\Delta.\]
We claim that the \'etale $\Zp$-local systems $\bL_n\coloneqq q_* ((\bA_{n,X}^{\Delta})^\vee )$ on $S$ detect the good reduction property of the family of curves $q:X\rightarrow S$. (In fact, it suffices to consider $n=5$.)

Note that the same construction above applies to the case where $S$ is a point. More precisely, for any geometrically connected smooth projective curve $Z$ over $K$, we obtain an \'etale $\Zp$-local system $\bA_{n,Z}^\Delta$ on $Z$. Consider
\[
C_n(Z)\coloneqq\mH^0_\et (Z_{\barK},(\bA_{n,Z}^\Delta)^\vee )[\frac 1p]
\]
viewed as a $\Gal_K$-representation. Applying to the situation $Z=X_s$ for any classical point $s$ of $S$, we have a base change isomorphism
\[
 (\bL_n)_{\bar s}[\frac 1p]\cong C_n(X_s)
\]
as $\Gal_{\kappa(s)}$-representations (cf. Proposition~\ref{prop: relative diagonal detector}). Therefore, by Theorem~\ref{thm: introduction crystalline rigidity} and Theorem~\ref{thm: introduction crystalline rigidity(semistable reduction)}, it reduces to showing that $C_n(Z)$'s detect the good reduction property of $Z$. 

Now, let $Z$ be a geometrically connected smooth projective curve over $K$ of genus $g\geq 2$, and we assume that $Z$ has semistable reduction over $\mathcal{O}_K$. Choose a closed point $b\in Z(K)$ and a geometric point $\bar{b}$ above it.\footnote{
Lemma~\ref{lem: descent along field extension} allows us to detect  the good reduction property of $Z$ after any finite base field extension. So we can assume that $Z$ has a $K$-point. Note that this step needs the condition $g\geq 2$.
} Consider the universal unipotent local systems $(\bE_n,u_n)$ in $\Uni_{\bZ_p}^{\leq n}(Z, b)$ as in \S \ref{subsection: intro_proof_of_main_thm_section}. In particular, we have an isomorphism
\begin{equation}\label{eq: into algebraic structure of local system}
 \mathrm{ev}_{u_n}\colon q_*\sEnd(\bE_n)\xra[]{\sim}b^*\bE_n
\end{equation}
by \eqref{eq: intro universality}. Let $A_n\coloneqq\bE_{n,\bar{b}}$, viewed as (a $\Zp$-lattice inside) a $\Gal_K$-representation. Then \eqref{eq: into algebraic structure of local system} induces an isomorphism $\mH_\et^0 (Z_{\barK},\sEnd(\bE_n) ) \cong A_n$ via the identification $\Loc_{\Zp}(\Spec K) \cong \Rep_{\Zp}(\Gal_K)$. This equips $A_n$ with a structure of $\bZ_p$-algebra. Consider its commutator quotient (i.e., the zero-th Hochschild homology)
\[
\HH_0(A_n)= A_n/[A_n,A_n] 
\]
where $[A_n, A_n]$ stands for the $\Zp$-submodule generated by the commutators. It carries a $\Gal_K$-action induced from that on $A_n$. It turns out $C_n(Z)$ is precisely computed by this commutator quotient.
\begin{prop}[Proposition~\ref{prop: diagonal vs adjoint}]
Let $Z$, $b$, and $\bar{b}$ be as above. There is an isomorphism 
\[
 C_n(Z)^\vee \cong \HH_0(A_n)[\frac 1p]
\] 
of $\Gal_K$-representations.
\end{prop}
Granting this, it remains to show the following:
\begin{thm}[{Theorem~\ref{thm: HH_0 good reduction}}]\label{thm: introduction HH_0 detects good reduction}
    Let $Z$, $b$, and $\bar{b}$ be as above. Then $Z$ has good reduction if and only if the $\Gal_K$-representation $\HH_0(A_5)[\frac 1p]$ is crystalline.
\end{thm}

In the rest of the introduction, we sketch the proof of Theorem~\ref{thm: introduction HH_0 detects good reduction}, using the machinery developed in \cite{CDPS}.
First, using formal deformation, we can construct a \emph{minimal semistable log curve} (cf. \S \ref{subsec: log de Rham fundamental group}) $\ul Y$ over the standard log point
\[
 \ul T= (\Spec K_0, \bN , 1\mapsto 0)^a
\]
(where $K_0=W(k)[\frac{1}{p}]$) with a good section $b_0$\footnote{A section is \emph{good} if it lies in the smooth locus of $Y$.}, such that $Z$ has good reduction if and only if the underlying curve $Y$ of $\ul{Y}$ is smooth. 

Secondly, the construction in \cite[\S2--\S 3]{CDPS} produces the \emph{relatively unipotent log de Rham fundamental group} $\cG_Y$ of $(\ul Y/\ul T,b_0)$ equipped with a monodromy operator $N$. Moreover, there is a relatively universal projective system of integrable log connections $\{(\cW_{n,Y}, w_{n,Y})\}_{n\geq 1}$ (cf. Proposition~\ref{prop: CDPS universal objects}) satisfying
\[
  \cG_Y \cong \sSpec (\varinjlim_n(b_0^*\cW_{n,Y})^\vee ),
\]
(also see \cite[Theorem 3.12]{CDPS}). Let $B_{n,Y}= b_0^*\cW_{n,Y}$ and denote its derivation by $D_n$. Moreover, let $\kg_Y$ be the Lie algebra of $\cG_Y$. Consider its universal enveloping algebra $U(\kg_Y)$ with two-sided ideal $I$ generated by $\kg_Y$. We prove in Proposition~\ref{prop: Lie algbebra vs B_n} that there is a natural isomorphism
\begin{equation}\label{eq: intro lie algebra}
    B_{n,Y}^{\mathrm{op}} \cong  U(\kg_Y)/I^n
\end{equation}
compatible with derivations.

As in the case of $A_n$, the module $B_{n,Y}$ is naturally a $K_0$-algebra so that we can consider its commutator quotient $\HH_0(B_{n,Y})$. In Proposition~\ref{prop: semistable comparison}, we show that there is a natural isomorphism
\[
D_{\st,K} (\HH_0(A_{n,Z})[\frac{1}{p}] ) \cong \HH_0(B_{n,Y}), 
\]
intertwining the monodromy on the left-hand side and the derivation on the right-hand side.
Therefore, we are left to show that $Y$ is smooth if and only if $D_5=0$ on $\HH_0(B_{5,Y})$. 

If $Y$ is smooth, then $N$ is trivial on $\cG_Y$ by \cite[Proposition 6.8]{CDPS}, hence $D_5=0$ on $\HH_0(B_{5,Y})$. Conversely, suppose $Y$ is not smooth, we want to argue that $D_5\neq 0$. Consider the dual graph $\Gamma$ of $Y$. There are two cases:
\begin{enumerate}
    \item If $\Gamma$ contains a loop, then $N$ is non-trivial on $\cG_Y^{\ab}$ by \cite[Theorem 6.24]{CDPS}, which implies that $D_2 \neq 0$ on $B_{2,Y}= \HH_0(B_{2,Y})$. Hence, $D_5$ is also non-trivial on $\HH_0(B_{5,Y})$.
    \item Otherwise, $\Gamma$ must be a tree. In this case, the pushout construction in \cite[Proposition 6.25]{CDPS} yields a minimal semistable curve $W$, whose dual graph is a chain
    \[
     W_0-W_1-\cdots -W_s,
    \]
    in which $W_i$ is isomorphic to $\bP^1$ for each $0< i< s$. Moreover, there is a natural surjection $\kg_{Y} \twoheadrightarrow \kg_{W}$ between Lie algebras (cf. Proposition~\ref{prop: contraction to chain}), which further induces a surjection $\HH_0(B_{5,Y}) \twoheadrightarrow \HH_0(B_{5,W})$ compatible with derivations. Using \eqref{eq: intro lie algebra} and \cite[Remark 6.28]{CDPS}, one can directly calculate that $D_5\neq 0$ on $\HH_0(B_{5,W})$ (cf. Lemma~\ref{lem: chain HH_0 trace}). It follows that $D_5 \neq 0$ on $\HH_0(B_{5,Y})$.
\end{enumerate}
This finishes the proof of Theorem~\ref{thm: introduction HH_0 detects good reduction} and hence the main theorems.

\vspace{0.1in}
\subsection*{Notation and conventions}
\begin{itemize}
\item Throughout the article, $K$ denotes a $p$-adic field, with ring of integers $\mathcal{O}_K$, uniformizer $\varpi$, and residue field $k$. Let $\mathbb{C}_K=\widehat{\overline{K}}$ and $K_0=W(k)[\frac 1p]$.
\item Let $\mathrm{Rep}_{\Qp}(\Gal_K)$ (resp. $\mathrm{Rep}_{\Zp}(\Gal_K)$) denote the category of continuous representations of $\Gal_K$ on finite-dimensional $\Qp$-vector spaces (resp. on finite free $\Zp$-modules).
\item Rigid analytic varieties over $K$ are viewed as adic spaces over $\mathrm{Spa}(K, \mathcal{O}_K)$. For a classical point $x$ on a rigid analytic variety, we use $\kappa(x)$ to denote its residue field. For an algebraic variety $X$ over $K$, let $X^{\mathrm{an}}$ denote its analytification. 
\item For a $p$-adic formal scheme $\mathfrak{X}$ over $\mathcal{O}_K$ that is topologically of finite type, we use $\mathfrak{X}_{\eta}$ to denote its rigid analytic generic fiber. We say $\mathfrak{X}$ is \emph{semistable} if it is strictly semistable in the sense that, \'etale locally, it admits an \'etale morphism
\[\mathfrak{X}\rightarrow \mathrm{Spf} (\mathcal{O}_K\langle T_0,\ldots, T_d\rangle/(T_0\cdots T_r-\varpi) )\]
for some integers $d\geq r\geq 0$. 
\item Similarly in the algebraic setting, we say that an algebraic scheme over $\mathcal{O}_K$ is \emph{semistable} if it is strictly semistable.
\item For any rigid analytic variety $X$ over $K$, we use the term \emph{\'etale $\Zp$-local system} for a torsion-free lisse $\Zp$-sheaf on $X_{\et}$ (cf. \cite[Definition 9.1.10]{kedlayaliu1}). This is the convention used in \cite{DY}, \cite{p_adic_monodromy_curve}, and \cite{relative}, but different from the one used in \cite{diao2019logarithmic}. Let $\mathrm{Loc}_{\Zp}(X)$ denote the category of \'etale $\Zp$-local systems on $X$. We use $\sHom$ and $\sEnd$ to denote the internal sheaves in this category, while we use $\Hom$, $\End$, and $\Ext^i$ for the corresponding global groups. Moreover, we use $\mathrm{Loc}_{\Qp}(X)$ to denote the category of \'etale $\Qp$-local systems on $X$.
\end{itemize}

\subsection*{Acknowledgment}
The author would like to express his sincere gratitude to his advisor, Hansheng Diao, for suggesting the topic of this thesis and for his guidance and support throughout the author’s doctoral studies. The author also thanks him warmly for many helpful discussions. He is grateful to Professor Diao for carefully reading earlier drafts of this thesis and offering valuable suggestions for improvement.

\vspace{0.2in}
\section{Universal unipotent local systems}
\label{sec: universal unipotent local systems}

The goal of this section is to construct a collection of universal unipotent $\Zp$-local systems associated with a family of curves; this generalizes the construction in \cite{AIK} to the relative setting. In this section, let $S$ be a geometrically connected smooth rigid analytic variety over $K$ and let $q:X\rightarrow S$ be a smooth proper family of geometrically connected curves. Assume that $q$ admits a section $b\colon S\rightarrow X$; in other words, $q$ is a family of \emph{pointed} curves.

\vspace{0.05in}
\subsection{Relatively unipotent local systems}\label{subsec: relative unipotent local systems}
\noindent
\vspace{0.05in}

\noindent

\begin{defn}\label{def: relative unipotence}
\begin{enumerate}
\item An \'etale $\Zp$-local system on $X$ is called \emph{relatively constant} if it is of the form $q^*\bV$ for an \'etale $\Zp$-local system $\bV$ on $S$.

\item For any integer $n\geq 1$, an \'etale $\Zp$-local system $\bL$ on $X$ is called \emph{relatively unipotent with unipotency index $\leq n$} if it admits a filtration
\[
 \bL=\bL^1\supseteq\bL^2\supseteq\cdots\supseteq\bL^n
 \supseteq\bL^{n+1}=0
\]
where each quotient $\bL^i/\bL^{i+1}$ is relatively constant. Let $\uni(\bL)$ denote the smallest possible length of such a filtration; in particular, $\uni(0)=0$. Let $\Uni_{\Zp}^{\leq n}(q)$ denote the full subcategory of $\mathrm{Loc}_{\Zp}(X)$ consisting of relatively unipotent ones with unipotency index $\leq n$, and put
\[
 \Uni_{\Zp}(q)\coloneqq \bigcup_{n\geq1}\Uni_{\Zp}^{\leq n}(q).
\]
\end{enumerate}
\end{defn}

Note that $\Uni_{\Zp}(q)$ is an exact tensor category. Taking duals preserves the unipotency index, and the tensor product filtration yields an inequality
\[
 \uni(\bL\otimes\bL')\leq\uni(\bL)+\uni(\bL')-1
\]
for $\bL,\bL'\in \Uni_{\Zp}(q)$.

\begin{defn}\label{def: pointed relative local systems}
For any integer $n\geq 1$, let $\Uni_{\Zp}^{\leq n}(q,b)$ denote the category consisting of pairs $(\bL,v)$ with $\bL\in\Uni_{\Zp}^{\leq n}(q)$ and
\[
 v \in \Hom(\ul\Zp,b^*\bL)= \mH^0_\et(S,b^*\bL).
\]
We refer to $v$ as a \emph{marking} of $\bL$. A morphism $(\bL_1,v_1)\ra(\bL_2,v_2)$ is a morphism $\phi:\bL_1\ra\bL_2$ such that $v_2=(b^*\phi)\circ v_1$. Then we put 
\[
 \Uni_{\Zp}(q,b)\coloneqq \bigcup_{n\geq1}\Uni_{\Zp}^{\leq n}(q,b).
\]
\end{defn}

\begin{defn}\label{def: relative universal system}
A projective system
\[
 \cdots\ra(\bE_{n+1},u_{n+1})\xra{\lambda_n}(\bE_n,u_n)
 \ra\cdots\xra{\lambda_1}(\bE_1,u_1)
\]
in $\Uni_{\Zp}(q,b)$ is called \emph{relatively universal} if $\uni(\bE_n)\leq n$ for every $n\geq 1$, and for any $\bL\in\Uni_{\Zp}^{\leq n}(q)$, evaluation at $u_n$ induces an isomorphism
\begin{equation}\label{eq: relative evaluation}
 \ev_{u_n}:q_*\sHom(\bE_n,\bL)\xra{\sim}b^*\bL.
\end{equation}
Here, the evaluation map $\ev_{u_n}$ is defined as the composition
\[
 q_*\sHom(\bE_n,\bL)
 \cong b^*q^*q_*\sHom(\bE_n,\bL)
 \ra b^*\sHom(\bE_n,\bL)\xra{u_n}b^*\bL,
\]
where the middle arrow is induced by the counit $q^*q_*\ra\id$.
\end{defn}

\begin{rmk}\label{rmk: pointed uniqueness}
Note that $(\bE_n,u_n)$ is an initial object in $\Uni_{\Zp}^{\leq n}(q,b)$. Indeed, for any $(\bL, v) \in \Uni_{\Zp}^{\leq n}(q,b)$, there is a unique morphism $(\bE_n ,u_n) \ra (\bL, v)$ given by \eqref{eq: relative evaluation} by taking global sections. In particular, the relative universal projective system, if exists, is unique up to a unique isomorphism. 
\end{rmk}

Following \cite[\S 3.5]{AIK}, we provide a way to check whether a projective system is relatively universal. Consider the following conditions on a projective system $\{(\bE_n,u_n, \lambda_n)\}_{n \geq 1}$:
\begin{enumerate}
 \item[(RU1)] $\bE_1=\ul\Zp$ and $u_1=1$.
 \item[(RU2)] For every $n\geq1$, the map $\lambda_n:\bE_{n+1}\rightarrow \bE_{n}$ fits into an exact sequence
 \begin{equation}\label{eq: universal extension shape}
 0\ra\bT_n\ra\bE_{n+1}\xra{\lambda_n}\bE_n\ra0,
 \end{equation}
 where $\bT_n=q^*\bB_n$ for some $\bB_n\in \mathrm{Loc}_{\Zp}(S)$.
 \item[(RU3)] The connecting morphism\footnote{Here, $\sExt^i_S$ stands for the relative Ext sheaf; namely, for any $\bF_1, \bF_2\in \mathrm{Loc}_{\Zp}(X)$, we put $\sExt^i_S(\bF_1, \bF_2)\coloneqq R^iq_*\sHom(\bF_1, \bF_2)$. When $i=0$, we put $\sHom_S(\bF_1, \bF_2)\coloneqq \sExt^0_S(\bF_1, \bF_2)=q_*\sHom(\bF_1, \bF_2)$.}
 \[
 \delta_n:\bB_n^\vee\cong\sHom_S(\bT_n,\ul\Zp)
 \ra\sExt_S^1(\bE_n,\ul\Zp)
 =\mR^1q_*\bE_n^\vee
 \]
 induced by \eqref{eq: universal extension shape} is an isomorphism.
\end{enumerate}

Since all fibers of $q$ are geometrically connected, we have $q_*\ul\Zp \cong \ul\Zp$. Dualizing \eqref{eq: universal extension shape} and using (RU3), we deduce that $q_*\bE_n^\vee\xra{\sim}q_*\bE_{n+1}^\vee$ and the map $\lambda_n^*:\mR^1q_*\bE_n^\vee\longrightarrow\mR^1q_*\bE_{n+1}^\vee$ is the zero map. In particular, the composition
\[
 \epsilon_n\coloneqq \lambda_1\circ\cdots\circ\lambda_{n-1}:\bE_n\ra\ul\Zp,
\]
induces an isomorphism 
\begin{equation}\label{eq: universal cohomology maps}
q_*\epsilon_n^\vee:\ul\Zp\xra{\sim}q_*\bE_n^\vee. 
\end{equation}
\begin{prop}\label{prop: universal axioms}
A projective system $\{(\bE_n,u_n, \lambda_n)\}_{n \geq 1}$ satisfying (RU1)--(RU3) is relatively universal.
\end{prop}

\begin{proof} We prove the following statement by induction on $n$:
For any integers $m\geq n\geq1$ and any $\bL\in\Uni_{\Zp}^{\leq n}(q)$, evaluation at $u_m$ induces an isomorphism
\[
 \ev_{u_m}:q_*\sHom(\bE_m,\bL)\xra{\sim}b^*\bL.
\]
When $n=1$, this follows from \eqref{eq: universal cohomology maps} and the projection formula. Now suppose $n>1$ and that the assertion holds at level $n-1$. For $\bL \in \Uni_{\Zp}^{\leq n}(q)$, choose a presentation
\begin{equation}\label{eq: presentation of L}
 0\ra q^*\bV\ra\bL\ra\bL'\ra0
\end{equation}
for some $\bL'\in\Uni_{\Zp}^{\leq n-1}(q)$ and $\bV\in \mathrm{Loc}_{\Zp}(S)$. Note that pulling back along $\lambda_{m-1}$ induces an isomorphism
\[
 \lambda_{m-1}^*\colon q_*\sHom(\bE_{m-1},\bL')\xra[]{\sim}
 q_*\sHom(\bE_{m},\bL').
\]
Consider the commutative diagram
\[
\begin{tikzcd}[column sep=large]
 q_*\sHom(\bE_{m-1},\bL')\ar[r,"\partial_{m-1}"]\ar[d,"\lambda_{m-1}^*"',"\cong"]
 &\sExt_S^1(\bE_{m-1},q^*\bV)\ar[d,"\lambda_{m-1}^*"]\\
 q_*\sHom(\bE_{m},\bL')\ar[r,"\partial_{m}"]
 &\sExt_S^1(\bE_{m},q^*\bV)
\end{tikzcd}
\]
where $\partial_{m-1}$ and $\partial_{m}$ are the connecting maps induced by \eqref{eq: presentation of L}. Since $\lambda_{m-1}^*\colon R^1q_*\bE^{\vee}_{m-1}\rightarrow R^1q_*\bE^{\vee}_m$ is the zero map, the right vertical arrow in the diagram is zero by the projection formula. This implies $\partial_{m}=0$. Consequently, evaluation at $u_m$ induces a diagram
\[
\begin{tikzcd}[column sep=small]
0\ar[r]&\bV\ar[r]\ar[d,equal]&q_*\sHom(\bE_{m},\bL)\ar[r]\ar[d,"\ev_{u_{m}}"]&q_*\sHom(\bE_{m},\bL')\ar[r]\ar[d,"\ev_{u_{m}}"]&0\\
0\ar[r]&\bV\ar[r]&b^*\bL\ar[r]&b^*\bL'\ar[r]&0
\end{tikzcd}
\]
with exact rows. By the induction hypothesis, the right vertical arrow is an isomorphism, and it follows that the middle vertical arrow is an isomorphism as well.
\end{proof}


\vspace{0.05in}
\subsection{Construction of the universal unipotent local systems}
\noindent
\vspace{0.05in}

\noindent

\begin{thm}\label{thm: existence of relative universal system}
There exists a relatively universal projective system $\{(\bE_n,u_n, \lambda_n)\}_{n\geq 1}$. Moreover, for each $n\geq 1$,  $\mR^1q_*\bE_n^\vee$ is an \'etale $\ul\Zp$-local system on $S$ (i.e., it is torsion-free).
\end{thm}

\begin{proof}
Following the strategy in \cite[\S 3.6]{AIK}, we construct a projective system satisfying (RU1)--(RU3) by induction on $n$. Since the fibers of $q$ are geometrically connected smooth projective curves, it follows that $q_*\ul\Zp\cong \ul\Zp$ and $\mR^1q_*\ul\Zp$ is torsion-free. We put $\bE_1=\ul\Zp$ and $u_1=1$. 

Suppose that the system has been constructed up to level $n$, satisfying (RU1)--(RU3) at all these levels, and that $\mR^1q_*\bE_n^\vee$ is torsion-free. Let $\bB_n\coloneqq (\mR^1q_*\bE_n^\vee)^\vee$ and $\bT_n\coloneqq q^*\bB_n$. We shall construct $\bE_{n+1}$ by computing the extension class. By the Leray spectral sequence (cf. \cite[Tag 0732]{stacks}), we have 
\[
 E_2^{i,j}=\mH^i_\et (S,\sExt_S^j(\bE_n,\bT_n) )
 \ \Rightarrow\ \Ext^{i+j}(\bE_n,\bT_n)
\]
which leads to a five-term exact sequence
\begin{equation}\label{eq: extension edge sequence}
\begin{split}
0\rightarrow\mH^1_{\et}(S, q_*\sHom(\bE_n, \bT_n))&\xrightarrow{\alpha_1}\mH^1_{\et} (X, \sHom(\bE_n, \bT_n) )\xrightarrow{\beta_1} \mH^0_\et (S,\sExt_S^1(\bE_n,\bT_n) )\\
&\rightarrow\mH^2_{\et} (S, q_*\sHom(\bE_n, \bT_n) )\xrightarrow{\alpha_2} \mH^2_{\et} (X, \sHom(\bE_n, \bT_n) )
\end{split}
\end{equation}
Using \eqref{eq: universal cohomology maps} and the projection formula, \eqref{eq: extension edge sequence} becomes
\[
0\ra\mH^1_\et(S,\bB_n)\xra{\alpha_1}\Ext^1(\bE_n,\bT_n)
 \xra{\beta_1}\End(\bB_n)\ra\mH^2_\et(S,\bB_n)\xra{\alpha_2}\Ext^2(\bE_n,\bT_n).
\]
Note that in the third term we have
\[
 \mH^0_\et (S,\sExt_S^1(\bE_n,\bT_n) )=\mH^0_\et (S, \mR^1q_*\sHom(\bE_n, \bT_n) )\cong \mH^0_\et (S, \mR^1q_*\bE_n^{\vee}\otimes \bB_n )
 \cong\mH^0_\et (S,\bB_n^\vee\otimes\bB_n )
 =\End(\bB_n).
\]
Also note that the map $\alpha_i$ can be identified with the composition
\[
 \mH^i_\et(S,\bB_n)=\Ext^i(\ul\Zp,\bB_n)\xrightarrow[]{q^*}\Ext^i(q^*\ul\Zp,q^*\bB_n)\xrightarrow[]{\epsilon_n^*}\Ext^i(\bE_n,q^*\bB_n)
\]
for $i=1,2$. On the other hand, the marking $u_n:\ul\Zp\ra b^*\bE_n$ induces a map
\[
 r_i\coloneqq u_n^*\circ b^*:\Ext^i(\bE_n,q^*\bB_n)
 \ra\Ext^i(\ul\Zp,\bB_n).
\]
Since $(b^*\epsilon_n)\circ u_n=\id$, we have $r_i\circ \alpha_i=\id$. Thus, $\alpha_2$ is injective, and so $\beta_1$ is surjective. We arrive at an isomorphism
\[
 \beta_1:\ker(r_1)\xra{\sim}\End(\bB_n).
\]
Let $\xi_n$ be the unique class in $\ker(r_1)$ mapping to $\id_{\bB_n}$, and we choose an extension
\begin{equation}\label{eq: define E_n+1}
 0\ra\bT_n\ra\bE_{n+1}\xra{\lambda_n}\bE_n\ra0.
\end{equation}
representing $\xi_n$. By construction, $\bE_{n+1}$ belongs to $\Uni_{\Zp}^{\leq n+1}(q)$.
Since the image of $\xi_n$ under $r_1$ is zero, we obtain a commutative diagram
\[
\begin{tikzcd}
0\ar[r] & b^*\bT_n\ar[r]\ar[d,equal]& b^*\bE_{n+1}\ar[r, "b^*\lambda_n"] & b^*\bE_n\ar[r]&0\\
0\ar[r] & \bB_n \ar[r] &\bP_n \ar[r]\ar[u]&\ul\Zp \ar[r]\ar[u, "u_n"]& 0
\end{tikzcd}
\]
where the bottom row is split exact. A splitting yields $u_{n+1}:\ul\Zp\ra \bP_n\rightarrow b^*\bE_{n+1}$ with $(b^*\lambda_n)\circ u_{n+1}=u_n$. By unwinding the definitions, the identity $\beta_1(\xi_n)=\id_{\bB_n}$ ensures that the connecting map
\[
 \delta_n:\bB_n^\vee\ra\mR^1q_*\bE_n^\vee
\]
in (RU3) is an isomorphism. We have verified (RU1)--(RU3).

It remains to check $\mR^1q_*\bE_{n+1}^\vee$ is torsion-free. Dualizing \eqref{eq: define E_n+1} yields an exact sequence
\begin{equation}\label{eq: dual exact sequence}
 0\ra\bE_n^\vee\ra\bE_{n+1}^\vee\ra q^*\bB_n^\vee\ra0.
\end{equation}
Consider the long exact sequence induced by applying $q_*$ to \eqref{eq: dual exact sequence}. Note that the composition $\bB_n^\vee\ra q_*q^*\bB_n^{\vee}\rightarrow\mR^1q_*\bE_n^\vee$, up to the usual sign for dual extensions, is nothing but $\delta_n$, and is therefore an isomorphism. We obtain an injection
\[
 \mR^1q_*\bE_{n+1}^\vee
 \hra\mR^1q_*q^*\bB_n^\vee
 \cong(\mR^1q_*\ul\Zp)\otimes\bB_n^\vee
\]
where the last isomorphism follows from the projection formula. Consequently, $\mR^1q_*\bE_{n+1}^\vee$ is torsion-free.
\end{proof}


For later use, let us introduce the notion of augmented $\ul\Zp$-algebra object in the category of \'etale $\Zp$-local systems.

\begin{defn}
    Let $Y$ be any rigid analytic variety over $K$. An \emph{augmented $\ul\Zp$-algebra} in $\mathrm{Loc}_{\Zp}(Y)$ is an associative algebra object $\bA$ in $\mathrm{Loc}_{\Zp}(Y)$ with multiplication $m\colon \bA\otimes \bA\rightarrow \bA$ and unit $u\colon \ul\Zp \ra \bA$, together with an augmentation $\ve\colon \bA \ra \ul\Zp$ that commutes with multiplication, such that $\ve\circ u=\id$. The notion of \emph{augmented $\ul\Qp$-algebra} in $\mathrm{Loc}_{\Qp}(Y)$ is defined in a similar way with $\Zp$ replaced by $\Qp$.
\end{defn}

\begin{construction}\label{constr: universal fiber algebra}
Let $\{(\bE_n,u_n, \lambda_n)\}_{n\geq 1}$ be the relative universal projective system constructed in Theorem~\ref{thm: existence of relative universal system}. For each $n\geq 1$, we show that $\bA_n=b^*\bE_n\in \mathrm{Loc}_{\Zp}(S)$ is naturally an augmented $\ul\Zp$-algebra. Taking $\bL= \bE_n$ in \eqref{eq: relative evaluation}, we obtain an isomorphism
\[
    \mathrm{ev}_{u_n}\colon q_*\sEnd(\bE_n)\xra{\sim}\bA_n
\]
which equips $\bA_n$ with an associated algebra structure with $u_n(1)=\mathrm{ev}_{u_n}(\id)$. Moreover, the augmentation $\epsilon\colon \bA_n\rightarrow \ul\Zp$ is given by $b^*\epsilon_n\colon \bA_n=b^*\bE_n\rightarrow b^*\ul\Zp=\ul\Zp$. Since $b^*\lambda_n$ is compatible with multiplication, so is the augmentation. 
\end{construction}

\vspace{0.05in}
\subsection{The base change property}\label{subsec: the base change property}
\noindent
\vspace{0.05in}

\noindent
We show that relatively universal projective systems are compatible with base changes on $S$. For simplicity, we shall write $\{(\bE_n, u_n)\}$ instead of $\{(\bE_n, u_n, \lambda_n)\}$ when the context is clear.

\begin{prop}\label{prop: base change universal system}
Let $g:T\ra S$ be a morphism of smooth rigid analytic varieties over $K$ and consider the pullback diagram
\begin{equation}\label{diagram: pullback}
\begin{tikzcd}[column sep=1.2cm, row sep=1.2cm]
 X_T=X\times_S T\ar[r,"g_X"]\ar[d,"q_T"]&X\ar[d,"q"']\\
 T \ar[u, "b_T", bend left= 40pt]\ar[r,"g"]& S\ar[u, "b"', bend right= 40pt]
\end{tikzcd}
\end{equation}
where $b_T$ is the section of the family $q_T:X_T\rightarrow T$ induced from $b$. Let $\{(\bE_n,u_n)\}_{n\geq 1}$ be a relatively universal projective system in $\Uni_{\Zp}(q,b)$. Then $\{(g_X^*\bE_n,g^*u_n)\}_{n\geq 1}$ is a relatively universal projective system in $\Uni_{\Zp}(q_T,b_T)$. 
\end{prop}

\begin{proof}
We need to check (RU1)--(RU3) for the system $\{(g_X^*\bE_n,g^*u_n)\}_{n\geq 1}$. Pulling back along $g_X$ clearly preserves (RU1) and (RU2). In particular, the exact sequence in (RU2) has the form
\[
 0 \ra q_T^*g^*\bB_n\ra g_X^*\bE_{n+1}\xra{g_X^*\lambda_n} g_X^*\bE_n\ra0.
\]
The connecting map
\[
 \delta_n^T \colon g^*\bB_n^\vee \ra R^1q_{T,*}g_X^*\bE_n^\vee
\]
in (RU3) can be identified with the composition
\[
(g^*\bB_n)^\vee\xra{\,g^*\delta_n\,}
 g^*\mR^1q_*\bE_n^\vee
 \xra{\sim}\mR^1q_{T,*}g_X^*\bE_n^\vee
\]
where the last isomorphism is given by the proper base change theorem. Therefore, (RU3) also holds. 
\end{proof}

In later sections, we shall apply the base change theorem in the special case where $T$ is a classical point of $S$.

\begin{cor}\label{cor: base change classical point}
Let $s\in S$ be a classical point and consider the pullback diagram
\begin{equation}\label{diagram: pullback to point}
\begin{tikzcd}[column sep=1.6cm, row sep=1.2cm]
X_s \ar[r,"\iota_s"] \ar[d,"q_s"]
& X \ar[d,"q"']\\
s \ar[u, bend left=40pt,
      "b_s"{pos=0.5, yshift=6pt}]
  \ar[r,"s"']
& S \ar[u, bend right=40pt,"b"']
\end{tikzcd}
\end{equation}
where $b_s$ is the section of $q_s\colon X_s\rightarrow s$ induced from $b$. Let $\{(\bE_{n,s},u_{n,s})\}_{n\geq 1}$ be a relatively universal projective system in $\Uni_{\Zp}(q_s,b_s)$. Then there is a unique isomorphism
\[
 (\bE_{n,s},u_{n,s})\xra{\sim}(\iota_s^*\bE_n,s^*u_n).
\]
compatible with the transition maps.
\end{cor}

\begin{proof}
The pullback system $\{(\iota_s^*\bE_n,s^*u_n)\}_{n\geq 1}$ is a relatively universal projective system in $\Uni_{\Zp}(q_s,b_s)$ by Proposition~\ref{prop: base change universal system}. Then the assertion follows immediately from Remark~\ref{rmk: pointed uniqueness}.
\end{proof}

\vspace{0.05in}
\subsection{Comparison with the universal projective system of Andreatta--Iovita--Kim}\label{subsec: comparison with AIK}
\noindent
\vspace{0.05in}

\noindent
Let $Z$ be a geometrically connected smooth projective (algebraic) curve over $K$ of genus $g\geq2$, and let $b\in Z(K)$ be a $K$-rational closed point. Passing to rigid analytification, we obtain $q\colon Z^{\an}\ra\Spa(K,\cO_K)$ with a section $b$. Let $\{(\bE_n, u_n)\}_{n\geq 1}$ be the relatively universal projective system in $\Uni_{\Zp}(q,b)$ constructed in Theorem~\ref{thm: existence of relative universal system}. We remark that everything established in \S \ref{subsec: relative unipotent local systems}--\ref{subsec: the base change property} has a $\Qp$-version; that is, we simply replace $\Zp$-local systems by isogeny $\Zp$-local systems and replace $\ul\Zp$ by $\ul\Qp$. By a slight abuse of notation, we denote the category of unipotent isogeny $\Zp$-local systems by $\Uni_{\Qp}(q,b)$. In particular, $\{(\bE_n[\frac{1}{p}], u_n[\frac{1}{p}])\}_{n\geq 1}$ is a relatively universal projective system in $\Uni_{\Qp}(q,b)$.

On the other hand, Andreatta, Iovita, and Kim construct a universal projective system of (algebraic) $\Qp$-local systems on $Z_{\et}$ in \cite{AIK}, associated with the pair $(Z,b)$. More precisely, they first construct such a system over $Z_{\overline{K}, \et}$ in \cite[Proposition 3.5]{AIK}, and then descend it to a projective system $\{(\bF_n, v_n)\}_{n\geq 1}$ over $Z_{\et}$ using the $\Gal_K$-action (cf. \cite[Corollary 3.7]{AIK}). Passing to analytification, we obtain a projective system $\{(\bF_n^{\an}, v_n^{\an})\}_{n\geq 1}$ of unipotent isogeny $\Zp$-local systems on $Z^{\an}_{\et}$.

\begin{prop}\label{prop: algebraic analytic universal comparison}
For each $n\geq 1$, there is a unique isomorphism
\[
 \eta_n \colon  (\bE_n[\frac{1}{p}],u_n[\frac{1}{p}] )\xra{\sim}(\bF_n^{\an},v_n^{\an})
\]
in $\Uni_{\Qp}(q,b)$, and these isomorphisms are compatible with the transition maps. In particular, $b^*\eta_n$ identifies $\bA_n[\frac 1p]$ with $\bF_{n,\bar{b}}$ as augmented $\Qp$-algebra objects in $\mathrm{Loc}_{\Qp}(\mathrm{Spa}(K, \mathcal{O}_K))\cong \mathrm{Loc}_{\Qp}(\mathrm{Spec}\,K)\cong \mathrm{Rep}_{\Qp}(\Gal_K)$.
\end{prop}

\begin{proof}
We first establish such isomorphisms over $\mathbb{C}_K=\widehat{\barK}$. Note that for any $\Qp$-local system $\bL$ on $Z_{\overline{K}, \et}$ and any $i\geq 0$, there is a canonical isomorphism
\begin{equation}\label{eq: geometric algebraic analytic cohomology}
 \mH^i_\et(Z_{\barK},\bL)
 \xra{\sim}
 \mH^i_\et(Z_{\bC_K}^{\an},\bL^{\an})
\end{equation}
(cf. \cite[Theorem 3.7.2]{Huber1996Etale}). Let $(q_{\bC_K}, b_{\bC_K})$ be the base change of $(q,b)$ to $\bC_K$. By \cite[Proposition 3.5]{AIK} and \eqref{eq: geometric algebraic analytic cohomology}, $\{(\bF_n^{\an} , v_n^{\an})\}_{n\geq 1}$ satisfies (RU1)--(RU3) when viewed as a projective system in $\Uni_{\Qp}(q_{\bC_K}, b_{\bC_K})$. On the other hand, $\{(\bE_n[\frac{1}{p}], u_n[\frac{1}{p}])\}_{n\geq 1}$ also satisfies (RU1)--(RU3), viewed as a projective system in $\Uni_{\Qp}(q_{\bC_K}, b_{\bC_K})$. Consequently, by (the $\Qp$-version of) Remark~\ref{rmk: pointed uniqueness}, there are unique isomorphisms
\[
 \eta_{n,\bC_K}\colon  (\bE_n[\frac{1}{p}],u_n[\frac{1}{p}] )\xra{\sim}
 (\bF_n^{\an},v_n^{\an})
\]
in $\Uni_{\Qp}(q_{\bC_K}, b_{\bC_K})$, compatible with transition maps. 

Now, by looking at the descent data given by the $\Gal_K$-actions on both systems, we see that $\eta_{n, \bC_K}$ uniquely descends to an isomorphism
\[
 \eta_n\colon  (\bE_n[\frac{1}{p}],u_n[\frac{1}{p}] )\xra{\sim}(\bF_n^{\an},v_n^{\an})
\]
in $\Uni_{\Qp}(q, b)$. These isomorphisms are necessarily compatible with transition maps by (the $\Qp$-version of) Remark~\ref{rmk: pointed uniqueness}.

\end{proof}

\vspace{0.2in}
\section{An upgraded criterion for good reduction of curves}

Throughout this section, let $Z$ be a  geometrically connected smooth projective curve over $K$ of genus $g\geq 2$. We use $q\colon Z \ra \Spec K$ to denote the structure morphism. The goal of this section is to provide a variant (actually, an upgrade) of the non-abelian criterion of Andreatta--Iovita--Kim (cf. Theorem~\ref{thm: introduction AIK}), without the datum of a fixed section. The main theorem is Theorem~\ref{thm: diagonal good reduction}, which will play a fundamental role in the proof of Theorem~\ref{thm: introduction rigidity}.

\vspace{0.05in}
\subsection{The diagonal construction}
\label{subsec: diagonal construction}

\begin{construction}\label{constr: diagonal universal system}
Consider morphisms
\[
 \pr_1:Z\times_K Z\ra Z,\ \
 \Delta:Z\ra Z\times_K Z,
\]
where $\pr_1$ stands for the projection to the first component, while $\Delta$ is the diagonal map. The pair $(\pr_1, \Delta)$ is a smooth projective family of geometrically connected pointed curves, as in the beginning of \S \ref{sec: universal unipotent local systems}. Applying Theorem~\ref{thm: existence of relative universal system} to  $(\pr_1, \Delta)$, we obtain a relatively universal projective system $\{(\bE_n^\Delta, u_n^\Delta)\}_{n\geq 1}$ of relatively unipotent local systems. Put
\begin{equation}\label{eq: diagonal local systems}
  \bA_n^\Delta= \Delta^*\bE_n^\Delta,
\end{equation}
then $\bA_n^\Delta$ is an augmented $\ul\Zp$-algebra in the category of \'etale $\bZ_p$-local systems on $Z$, with multiplication $\mu_n^\Delta$ and augmentation $\ve_n: \bA_n^\Delta \ra \ul\Zp$ (cf. Construction~\ref{constr: universal fiber algebra}). 
By construction, we have $\ve_n \circ u_n^\Delta=1$. 
Consider
\begin{equation}\label{eq: intrinsic arithmetic detectors}
 C_n(Z)\coloneqq \mH^0_\et(Z_{\barK},(\bA_n^\Delta)^\vee)[\frac 1p],
\end{equation}
viewed as an element in $\mathrm{Rep}_{\Qp}(\Gal_K)$.
\end{construction}

\begin{thm}
\label{thm: diagonal good reduction}
Let $Z$ be a geometrically connected  smooth projective curve over $K$ of genus $g \geq 2$. Suppose that $Z$ has semistable reduction over $\mathcal{O}_K$. Then $Z$ has good reduction over $\mathcal{O}_K$ if and only if $C_5(Z)$ is crystalline.
\end{thm}

The remainder of this section is dedicated to the proof of Theorem~\ref{thm: diagonal good reduction}. We first show that if $Z$ has good reduction (resp. semistable reduction), then $C_n(Z)$ is crystalline (resp. semistable) for all $n\geq 1$.

\begin{lem}\label{lem: pullback diagonal along section}
Suppose $b\in Z(K)$ is a closed point and consider the morphism
\[
i_b \colon Z \ra Z\times_K Z,
\]
sending $z \mapsto (b,z)$. Then the projective system $\{(i_b^*\bE_n^\Delta, b^*u_n^\Delta)\}_{n \geq 1}$ is relatively universal in $\Uni_{\Zp}(Z,b)$.
\end{lem}

\begin{proof}
This follows immediately by applying Proposition~\ref{prop: base change universal system} to the base change diagram
\begin{equation}\label{diagram: diagonal pullback}
\begin{tikzcd}[column sep=1.6cm, row sep=1.2cm]
 Z\ar[r,"i_b"]\ar[d,"q"]& Z\times_K Z\ar[d,"\pr_1"']\\
 \Spec K \ar[u, bend left=40pt, "b"]
  \ar[r, "b"]& Z.\ar[u, bend right=40pt, "\Delta"']
\end{tikzcd}
\end{equation}
\end{proof}

\begin{cor}\label{cor: good reduction implies crystalline}
    Suppose that $Z$ has good (resp. semistable) reduction over $\mathcal{O}_K$. Then $C_n(Z)$ is crystalline (resp. semistable) for every $n\geq 1$. 
\end{cor}

\begin{proof}
    Let $z$ be a closed point in $Z$ with residue field $L= \kappa(z)$. Consider the base change $Z_L$ and choose a closed point $\tilde{z}\in Z_L(L)$ that maps to $z$. Let $\overline z$ be a geometric point above $\tilde{z}$, which can also be viewed as a geometric point above $z$. Consider the pointed curve $(Z_L,\tilde{z})$ and let $\{(\bE_n, u_n)\}_{n\geq 1}$ be a relatively universal projective system in $\Uni_{\Zp}(Z_L,\tilde{z})$. By Lemma~\ref{lem: pullback diagonal along section} and Remark~\ref{rmk: pointed uniqueness}, there are unique isomorphisms
    \[(i_{\tilde{z}}^*\bE_n^\Delta, \tilde{z}^*u_n^\Delta)\cong (\bE_n, u_n)\]
    compatible with transition maps. Here, we have slightly abused the notation and view $\bE_n^\Delta$ as a local system on $Z_L\times_L Z_L$ via pullback. This implies
    \[
     \tilde{z}^*\bA_n^\Delta \cong \tilde{z}^*i_{\tilde{z}}^*\bE_n^\Delta \cong \tilde{z}^*\bE_n.
    \]
    Thus, as $\Gal_L$-representations, we have $\bA_{n, \overline z}^\Delta \cong \bE_{n, \overline z}$. Since $Z$ has good (resp. semistable) reduction over $\mathcal{O}_K$, $Z_L$ also has good (resp. semistable) reduction over $\mathcal{O}_L$ (see for example, \cite[Tag 0CDE]{stacks}). Therefore, $\bE_{n, \overline z}$ is crystalline (resp. semistable) by \cite[Theorem 1.4, Theorem 1.9]{AIK} and Proposition~\ref{prop: algebraic analytic universal comparison}. The same therefore holds for $\bA_{n, \overline z}^\Delta$. Consequently, by the pointwise criterion \cite[Theorem 7.2]{Guoyang}, the local system $\bA_n^\Delta$, as well as its dual $(\bA_n^\Delta)^\vee$, is crystalline (resp. semistable). 
    In the case of good reduction, by the crystalline comparison theorem \cite[Corollary C]{GuoReinecke2024}, the cohomology 
    \[
    C_n(Z)=\mH^0_\et (Z_{\barK},(\bA_n^\Delta)^\vee )[\frac 1p]
    \]
    of $(\bA_n^\Delta)^\vee$ is crystalline. In the case of semistable reduction, by the semistable comparison theorem \cite[Theorem 1.1]{log2} and \cite[Theorem 1.1]{DMS}, $C_n(Z)$ is semistable.
\end{proof}

\vspace{0.05in}
\subsection{Commutator quotients and adjoint local systems}
\noindent
\vspace{0.05in}

\noindent 
To prove the converse direction of Theorem~\ref{thm: diagonal good reduction}, we will describe $C_n(Z)$ in terms of a certain commutator quotient.

Recall that for any unital associative algebra $B$ over a commutative ring $R$, its zero-th Hochschild homology is defined to be the commutator quotient
\[
\HH_0(B) \coloneqq B/[B,B] = \mathrm{coker} (B\otimes_R B \xra{a\otimes c\, \mapsto \,ac- ca} B),
\]
where $[B,B]$ stands for the $R$-linear span of the commutators. 


\begin{lem}\label{lem: algebraic centralizer duality}
Let $A$ be a unital associative algebra over a commutative ring $R$, and let $E$ be a nonzero finite left $A$-module. Put $W=A^\vee\otimes_A E$, where duals are taken over $R$. Endow $A^\vee$ with the bi-$A$-module structure
\begin{equation}\label{eq: dual bimodule}
 (a\lambda a')(x)=\lambda(a'xa),
 \ \ \forall a,a',x\in A,\ \lambda\in A^\vee,
\end{equation}
and endow $W$ with the left $A$-action induced by the left action on $A^\vee$. Then,
\begin{enumerate}
 \item There are canonical morphisms of $A$-modules
 \begin{equation}\label{eq: algebraic middle maps}
 \begin{aligned}
    &\Hom_A(E,A^\vee)\otimes_A E \ra (\End_A(E))^\vee, \\ &\Hom_A(E,A^\vee)\otimes_A E \ra \Hom_A(E,W).  
 \end{aligned}
 \end{equation}
 When $E$ is a finite free $A$-module, they are both isomorphisms, which induce an isomorphism
 \begin{equation}\label{eq: algebraic centralizer duality}
     (\End_A(E))^\vee \cong \Hom_A(E,W).
 \end{equation}
 \item Suppose $E$ is finite free. For $\lambda\in A^\vee$, consider the $R$-linear map
 \[
  f_\lambda :E\ra W,\ \  v\longmapsto\lambda\otimes v.
 \]
  Then $f_\lambda$ is $A$-linear if and only if $\lambda \in \HH_0(A)^\vee$.
\end{enumerate}
\end{lem}

\begin{proof}
For (1), if we endow $E^\vee$ with the right $A$-module structure given by 
\[
 \lambda a (x)= \lambda(ax), \ \ \forall a\in A, \lambda \in E^\vee, x\in E,
\]
then there are natural morphisms
\begin{equation*}
\begin{aligned}
 & \alpha \colon E^\vee \otimes_A E \ra (\End_A(E))^\vee, \ \ &\lambda\otimes x \mapsto (T \mapsto \lambda(Tx)), \\
& \beta\colon \Hom_A(E,A^\vee) \ra E^\vee, \ \ &\lambda \mapsto (x\mapsto \lambda(x)(1)), \\
& \gamma \colon \Hom_A(E,A^\vee)\otimes_A E \ra \Hom_A(E,W), \ \ &T \otimes w \mapsto (x \mapsto T(x)\otimes w) .  
\end{aligned}
\end{equation*}
The desired morphisms \eqref{eq: algebraic middle maps} are given by $\alpha\circ  (\beta\otimes \id_E)$ and $\gamma$, respectively. When $E$ is free, $\alpha, \beta, \gamma$ are all isomorphisms.

For (2), since
\[
 f_\lambda(av)-a f_\lambda(v)
 =(\lambda a-a\lambda)\otimes v, \ \ \forall a\in A, v\in E,
\]
$f_\lambda$ is $A$-linear if and only if $\lambda a=a\lambda$ for every $a\in A$, or equivalently, $\lambda$ annihilates $[A,A]$.
\end{proof}

To proceed, we also need the notion of adjoint local systems.

\begin{defn}
For an \'etale $\Zp$-local system $\bL$ on $Z$, the adjunction $q^*q_*\sEnd(\bL)\ra\sEnd(\bL)$ makes $\bL$ a left $q^*q_*\sEnd(\bL)$-module. Define 
\[
 \bL^{\ad}\coloneq\sEnd_{q^*q_*\sEnd(\bL)}(\bL)
 \subseteq\sEnd(\bL),
\]
which is referred to as the \emph{adjoint local system} of $\bL$. It is naturally a $\underline{\Zp}$-algebra object in $\Loc_{\bZ_p}(Z)$.
\end{defn}

Now, suppose $Z$ contains a $K$-rational closed point $b\in Z(K)$ (this is the case when $Z$ has semistable reduction), and consider the pointed curve $(Z,b)$. Let $\{(\bE_n, u_n)\}_{n\geq 1}$ be a relatively universal projective system in $\Uni_{\Zp}(Z,b)$. Put $\bA_n= b^*\bE_n$, whose corresponding representation $\bE_{n,\bar{b}}$ in $\Rep_{\Zp}(\Gal_K)$ is denoted by $A_n$. The evaluation at $u_n$ induces 
\[
 \ev_{u_n}\colon q_*\sEnd(\bE_n) \xra[]{\sim} \bA_n.
\]
Thus, we have $\bE_n^\ad \cong \sEnd_{q^*\bA_n}(\bE_n)$.

\begin{prop}\label{prop: diagonal vs adjoint}
Let $Z$ and $b$ be as above. For any $n\geq 1$,
there is a canonical isomorphism of augmented $\ul{\bZ_p}$-algebras in $\Loc_{\Zp}(Z)$
\begin{equation}\label{eq: diagonal adjoint identification}
 \vartheta_{n,b}\colon \bA_n^\Delta= \Delta^*\bE_n^\Delta\xra{\sim}(\bE_n^\ad)^{\mathrm{op}},
\end{equation}
where $(\bE_n^\ad)^{\mathrm{op}}$ stands for the opposite algebra of $\bE_n^\ad$. Consequently, there is an isomorphism
\[
 C_n(Z)^\vee\cong\HH_0(A_n)[\frac 1p]
\]
as $\Gal_K$-representations.
\end{prop}

\begin{proof}
Put $r=q\circ \pr_1=q\circ \pr_2$ and define
\[
 \bM_n \coloneqq \sHom_{r^*\bA_n}(\pr_2^*\bE_n,\pr_1^*\bE_n).
\]
Pulling back along the diagonal yields
\[
 \Delta^*\bM_n=\sEnd_{q^*\bA_n}(\bE_n) =\bE_n^\ad.
\]
Taking $w_n= \id_{\bE_n}$ in $\mH_\et^0(Z,\Delta^*\bM_n)$ gives rise to an object $(\bM_n, w_n) \in \Uni_{\bZ_p}(\pr_1,\Delta)$. We claim that $\bM_n \in \Uni_{\Zp}^{\leq n}(\pr_1)$. 
Indeed, $\bM_n$ admits a decreasing filtration 
\[
\bM_n^i= \sHom_{r^*\bA_n}(\pr_2^*\bE_{n-i+1},\pr_1^*\bE_n), \ \ 1\leq i \leq n+1.
\]
For any $m\geq 1$, put $\bB_m\coloneqq (\mR^1q_*\bE_m^\vee)^\vee$. By \eqref{eq: universal extension shape},
the graded pieces of $\{\bM_n^i\}_{1\leq i \leq n+1}$ are given by  
\begin{align*}
  \bM_n^i/\bM_n^{i+1} &\cong \sHom_{r^*\bA_n}(\pr_2^*q^*\bB_{n-i}, \pr_1^*\bE_n) \\
  & \cong \sHom_{r^*\bA_n}(\pr_1^*q^*\bB_{n-i}, \pr_1^*\bE_n)
  \cong   \pr_1^*\sHom_{q^*\bA_n}(q^*\bB_{n-i},\bE_n) \
\end{align*}
which are relatively constant. This proves the claim.

The universality of $(\bE_n^\Delta,u_n^\Delta)$ induces a unique morphism $\phi_n:\bE_n^\Delta\ra\bM_n$ satisfying $\Delta^*\phi_n(u_n^\Delta)=w_n$.
Let $i_b:Z\ra Z\times_K Z$ be the closed embedding $z\mapsto(b,z)$ as above. Evaluation at $u_n \in \Hom(\ul\Zp,\bA_n)$ gives an isomorphism
\begin{equation}\label{eq: restrict to a point}
   i_b^*\bM_n
 \cong \sHom_{q^*\bA_n}(q^*b^*\bE_n,\bE_n) \cong \sHom_{q^*\bA_n}(q^*\bA_n,\bE_n)
 \xra[\ev_{u_n}]{\sim}\bE_n
\end{equation}
carrying the marking $b^*w_n= \id_{\bA_n}$ to $u_n$. By uniqueness, the composition of $i_b^*\phi_n$ and \eqref{eq: restrict to a point} coincides with the canonical isomorphism given in Lemma~\ref{lem: pullback diagonal along section}. In particular, $i_b^*\phi_n$ must be an isomorphism. Since $Z\times_K Z$ is connected, $\phi_n$ is also an isomorphism. Then
\[
    \vartheta_{n,b}\coloneqq \Delta^*\phi_n \colon \bA_n^\Delta\xra{\sim}\bE_n^\ad
\]
provides the desired isomorphism. We still have to check that the isomorphism preserves multiplications. 
Let $\mu^\ad$ denote the multiplication on $\bE_n^\ad$ induced by the composition. Let $\tau$ be the endomorphism of $\bA_n^\Delta\otimes_{\ul \Zp}\bA_n^\Delta$ sending $a\otimes c \mapsto c \otimes a$. We need to check that 
\begin{equation}\label{eq: one point linear}
    \psi\coloneqq  (\vartheta_{n,b}\circ \mu_n^\Delta\circ \tau - \mu^{\ad} \circ(\vartheta_{n,b} \otimes \vartheta_{n,b}) )\colon \bA_n^\Delta \otimes_{\ul\Zp} \bA_n^\Delta \ra \bE_n^{\ad}
\end{equation}
is a zero map. Since $Z$ is connected, we only need to check that $\psi_{\bar{b}}=0$. Consequently, we are left to check the isomorphism
\[
\vartheta_{\bar{b}} \coloneqq \vartheta_{n,b,\bar{b}} \colon A_n \xra[]{\sim} \End_{A_n}(A_n)
\]
is compatible with multiplication. Let $1_n\coloneqq u_n(1)\in A_n$. For any $a\in A_n$, we have $\vartheta_{\bar{b}}(a)(1_n)= a$ by construction. As $\vartheta_{\bar{b}}(a) \in \End_{A_n}(A_n)$, this implies
\[
 \vartheta_{\bar{b}}(a)(v)=  \vartheta_{\bar{b}}(a)(v\cdot 1_n)=v\cdot  \vartheta_{\bar{b}}(a)(1_n) =v\cdot a, 
\]
for all $v\in A_n$. Therefore, $\vartheta_{\bar{b}}(a) \vartheta_{\bar{b}}(c)= R_a \circ R_c= R_{ca}= \vartheta_{\bar{b}}(ca)$, where $R$ stands for the right multiplication.
Consequently, $\vartheta_{n,b}$ induces an isomorphism between $\bA_n^\Delta$ and $(\bE_n^\ad)^{\mathrm{op}}$, as augmented $\ul\Zp$-algebras.

It remains to show that $C_n(Z)^\vee \cong \HH_0(A_n[\frac 1p])$. Write $\bW= q^*\bA_n^\vee \otimes_{q^*\bA_n} \bE_n$, where $q^*\bA_n^\vee$ is equipped with a bimodule structure as in \eqref{eq: dual bimodule}. By \eqref{eq: algebraic middle maps} in Lemma~\ref{lem: algebraic centralizer duality} (1), there are natural morphisms
\[
 \begin{aligned}
    & \alpha\colon \sHom_{q*\bA_n}(\bE_n, q^*\bA_n^\vee)\otimes_{q^*\bA_n} \bE_n\ra (\bE_n^\ad)^\vee, \\ & \beta \colon \sHom_{q*\bA_n}(\bE_n,q^*\bA_n^\vee)\otimes_{q^*\bA_n} \bE_n \ra \sHom_{q^*\bA_n}(\bE_n, \bW).  
 \end{aligned}
\]
They become isomorphisms after pulling back along $b$ by \eqref{eq: algebraic centralizer duality} in Lemma~\ref{lem: algebraic centralizer duality} (1). Hence, both $\alpha $ and $\beta$ are isomorphisms. We arrive at a natural isomorphism
\[
  (\bE_n^\ad)^\vee \cong \sHom_{q^*\bA_n}(\bE_n, \bW)
\]
which further induces a natural inclusion
\begin{align*}
q_*(\bA_n^{\Delta})^\vee \cong q_*(\bE_n^\ad)^\vee \cong q_* (\sHom_{q^*\bA_n}(\bE_n, \bW) )\hookrightarrow q_*\sHom (\bE_n, \bW )
\xra[\ev_{u_n}]{\sim} b^*\bW \cong \bA_n^\vee
\end{align*}
between $\Zp$-local systems on $(\Spec K)_\et$. This induces 
\[
  \mH_\et^0 (Z_{\barK}, (\bA_n^{\Delta})^\vee ) \hookrightarrow A_n^\vee.
\]For any $\lambda\in A_n^\vee = \bA_{n,\bar{b}}^\vee$, the inverse of $\ev_{u_n}$ sends $\lambda$ to the morphism 
\[
f_\lambda \coloneqq \lambda \otimes\id_{\bE_n}  \colon \bE_n|_{Z_\barK} \ra \bW|_{Z_\barK}.
\]
in $\mH_\et^0 (Z_\barK,\sHom(\bE_n, \bW) )$. We need to determine when $f_\lambda$ belongs to $\mH_\et^0 (Z_\barK,\sHom_{q^*\bA_n}(\bE_n, \bW) )$. Since $Z$ is geometrically connected, we only need to analyze the $A_n$-linearity at $\bar b$. Note that $\bE_{n,\bar b}$ is a free $A_n$-module of rank one. According to Lemma~\ref{lem: algebraic centralizer duality} (2), $f_\lambda$ is $A_n$-linear if and only if $\lambda \in \HH_0(A_n)^\vee$. This finishes the proof by taking duals and inverting $p$.
\end{proof}

\begin{rmk}\label{rmk: geometric End in criterion}
For two closed points $b_1,b_2\in Z(K)$, write $\{(\bE_n^{b_1},u_n^{b_1})\}_{n\geq 1}$ and $\{(\bE_n^{b_2},u_n^{b_2})\}_{n\geq 1}$ for the corresponding relatively universal projective systems on $Z$. The identifications \eqref{eq: diagonal adjoint identification} induce a canonical change-of-base-point isomorphism
\[
 \vartheta_{n,b_2}\circ\vartheta_{n,b_1}^{-1}:
 \bE_n^{b_1,\ad}\xra{\sim}\bE_n^{b_2,\ad}.
\]
\end{rmk}

\begin{lem}\label{lem: descent along field extension}
    Let $Z$ be a geometrically connected smooth projective curve over $K$ of genus $g\geq 2$. Assume that $Z$ has semistable reduction over $K$. If $Z$ potentially has good reduction (i.e., $Z_L$ has good reduction over $\mathcal{O}_L$ for some finite extension $L/K$), then $Z$ has good reduction over $\mathcal{O}_K$.
\end{lem}
\begin{proof}
    Since $g\geq 2$, by \cite[Tag 0E8D]{stacks}, there are stable models $\cZ^{(K)}$ of $Z$ over $\cO_K$ and $\cZ^{(L)}$ of $Z_{L}$ over $\cO_L$. By the uniqueness of the stable model (cf. \cite[Tag 0E97]{stacks}),  there exists a unique isomorphism
    \[
     \cZ^{(K)}\times_{\Spec \cO_K}\Spec \cO_L \xra{\sim} \cZ^{(L)}.
     \]
     Since $Z_L$ has good reduction over $\mathcal{O}_L$, $\cZ^{(L)}$ is smooth. It follows that $\cZ^{(K)}$ is also smooth by faithful flat descent.
\end{proof}

Given these preparations, the proof of Theorem~\ref{thm: diagonal good reduction} reduces to the following criterion for good reduction. Its proof will be given in the next section.

\begin{thm}\label{thm: HH_0 good reduction}
Let $Z$ be a geometrically connected smooth projective curve over $K$ of genus $g\geq 2$. Suppose $Z$ has semistable reduction over $\mathcal{O}_K$. Suppose $Z$ contains a $K$-rational closed point $b\in Z(K)$ and let $\bar{b}$ be a geometric point above it. Let $\{(\bE_n,u_n)\}_{n \geq 1}$ be a relative universal projective system in $\Uni_{\Zp}(Z,b)$ and let $A_n= \bE_{n,\bar{b}}$. Then $Z$ has good reduction over $\mathcal{O}_K$ if and only if $\HH_0(A_5)[\frac 1p]$ is crystalline.
\end{thm}

\begin{proof}[Proof of Theorem~\ref{thm: diagonal good reduction} assuming Theorem~\ref{thm: HH_0 good reduction}]
By Corollary~\ref{cor: good reduction implies crystalline}, we only need to check the ``if'' part. Pick any closed point $b\in Z$ with residue field $L=\kappa(b)$. By Proposition~\ref{prop: diagonal vs adjoint} and Theorem~\ref{thm: HH_0 good reduction}, we know that $Z_L$ has good reduction over $\mathcal{O}_L$. Consequently, $Z$ also has good reduction over $\mathcal{O}_K$ by Lemma~\ref{lem: descent along field extension}.
\end{proof}

\vspace{0.2in}
\section{Good reduction in terms of commutator quotient}
The goal of this section is to prove Theorem~\ref{thm: HH_0 good reduction}, thereby completing the proof of Theorem~\ref{thm: diagonal good reduction}. Throughout the section, all curves are assumed to have genus at least two.

\vspace{0.05in}
\subsection{Relative unipotent log de Rham fundamental group}
\label{subsec: log de Rham fundamental group}
\noindent
\vspace{0.05in}

\noindent 
A key ingredient in the proof of Theorem~\ref{thm: HH_0 good reduction} is the theory of log de Rham fundamental groups developed in \cite{CDPS}. Throughout \S \ref{subsec: log de Rham fundamental group}, let $F$ be a field of characteristic zero and fix an algebraic closure $\overline F$. For any fs log scheme $\ul W$ over $F$, let $\MIC^{\nil}(\ul W/F)$ denote the category of locally free $\cO_W$-modules of finite rank with nilpotent integrable log connection. 
For simplicity, an object $(\cE,\nabla_\cE) \in \MIC^{\nil}(\ul W/F)$ will also be denoted by $\cE$ when the context is clear.
Let 
\[
 \ul T= (T,\cM_T) \coloneqq (\Spec F,\bN, 1\mapsto 0)^a
\]
be the standard log point over $F$, where $(-)^a$ stands for the log structure associated with the pre-log structure. Let $t\in \cM_T$ be the image of $1\in \bN$. Then we have $\Omega_{\ul T/ F}^1= F \cdot d\log t$. Let $\mathbf1_{T/F}\coloneqq (\cO_T,d)$ denote the trivial log connection on $T$. The category $\MIC^{\nil}(\ul T/F)$ is equivalent to the category of pairs $(V,N_V)$, where $V$ is a finite-dimensional vector space over $F$ and $N_V\in \End_F(V)$ is a nilpotent endomorphism .

We consider a geometrically connected projective log curve
\[
 q_0:\ul Y=(Y,\cM_Y)\ra\ul T.
\]
Here, by a \emph{log curve} we mean that $q_0$ a log smooth integral morphism of fs log schemes such that $Y$ is geometrically a reduced curve, as in \cite[Definition 6.1]{CDPS}. In particular, $\ul Y$ has only a vertical log structure.
We say $\ul Y$ is \emph{semistable} if \'etale locally, it admits an \'etale morphism
\[
Y\rightarrow \Spec(F[u,v]/(uv), \bN^2)^a,
\]
where the pre-log structure $\bN^2$ on $F[u,v]/(uv)$ is given by sending $(1,0)\mapsto u$ and $(0,1)\mapsto v$.

\begin{defn}[{cf. \cite[Definition 6.11]{CDPS}}]\label{def: minimal semistable}
    A semistable log curve $\ul Y$ as above is \emph{minimal semistable} if it is not a smooth curve of genus $1$, and every irreducible component of $Y_{\overline F}$ that is birationally equivalent to $\bP_{\overline F}^1$ must intersect with other components at more than one point.
\end{defn}

Moreover, a log curve $q_0\colon \ul Y \ra \ul T$  together with a section $b_0 \colon \ul T\ra\ul Y$ is referred to as a \emph{pointed log curve}. The section is called a \emph{good section} if $b_0$ lies in the smooth locus of $Y$.

Following \cite[\S 1]{CDPS}, we introduce the notion of relatively unipotent connections.

\begin{defn}\label{def: relative de Rham categories}
Let $(q_0 \colon \ul Y \ra \ul T, b_0)$ be a pointed log curve. 
\begin{enumerate}
\item A locally free module with integrable log connection $(\cE, \nabla_\cE)$ on $\ul Y$ is called \emph{relatively constant} if it is of the form $q_0^*(\cV,\nabla_{\cV})$ where $(\cV, \nabla_\cV)$ is a locally free module with integrable log connection on $\ul T$.

\item For any integer $n\geq 1$, a locally free module with integrable log connection $(\cE, \nabla_\cE)$ on $\ul Y$ is called \emph{relatively unipotent with unipotency index $\leq n$} if it admits a filtration
\[
 \cE= \cE^1\supseteq\ \cE^2 \supseteq  \cdots\supseteq \cE^n \supseteq \cE^{n+1}=0,
\]
where each quotient $\cE^i/\cE^{i+1}$ is relatively constant. Let $\uni(\cE,\nabla_{\cE})$ denote the smallest possible length of such a filtration; in particular, $\uni(0)=0$. Let $\Uni_{\dR}^{\leq n}(q_0)$ denote the full subcategory of $\MIC^{\nil}(Y/F)$ consisting of relatively unipotent ones with unipotency index $\leq n$, and put
\[
 \Uni_{\dR}(q_0)\coloneqq \bigcup_{n\geq1}\Uni_{\dR}^{\leq n}(q_0).
\]
\end{enumerate}
\end{defn}

\begin{defn}
For any integer $n\geq 1$, let $\Uni_{\dR}^{\leq n}(q_0,b_0)$ denote the category consisting of pairs $((\cE,\nabla_\cE),v)$ with $(\cE,\nabla_\cE)\in\Uni_{\dR}^{\leq n}(q_0)$ and
\[
 v \in \Hom_{\MIC^{\nil}(\ul T/F)} (\mathbf{1}_{T/F}, b_0^*(\cE,\nabla_\cE) ).
\]
A morphism $((\cE_1,\nabla_{\cE_1}),v_1) \ra ((\cE_2,\nabla_{\cE_2}),v_2)$ is a morphism $\phi:(\cE_1,\nabla_{\cE_1}) \ra (\cE_2,\nabla_{\cE_2})$ such that $(b_0^*\phi)(v_1) =v_2$. The section $v$ is called a \emph{marking} of $(\cE,\nabla_{\cE})$.
Then we put 
\[
 \Uni_{\dR}(q_0,b_0)\coloneqq \bigcup_{n\geq1}\Uni_{\dR}^{\leq n}(q_0,b_0).
\]
\end{defn}

Analogous to Definition~\ref{def: relative universal system}, there is a notion of relatively universal projective system in $\Uni_{\dR}(q_0,b_0)$.

\begin{defn}\label{def: relative universal connection}
A projective system
\[
 \cdots\ra(\cW_{n+1},w_{n+1})\xra{\lambda_n}(\cW_n,w_n)
 \ra\cdots\xra{\lambda_1}(\cW_1,w_1)
\]
in $\Uni_{\dR}(q_0,b_0)$ is called \emph{relatively universal} if $\uni(\cW_n)\leq n$ for every $n\geq 1$, and for any $\cE\in\Uni_{\dR}^{\leq n}(q_0)$, the evaluation at $w_n$ induces an isomorphism
\begin{equation}\label{eq: relative evaluation 2}
 \ev_{w_n} \colon q_{0,\dR,*}\sHom(\cW_n, \cE) \xra{\sim} b_0^* \cE
\end{equation}
in $\MIC^{\nil}(\ul T/ F)$. Here, $q_{0,\dR,*}\cF= \mH_\dR^0(\ul Y/ \ul T, \cF)$ denotes the relative horizontal sections, equipped with an induced connection on $\ul T$.
\end{defn}

\begin{rmk}\label{rmk: pointed uniqueness dR}
Following the same argument as in Remark~\ref{rmk: pointed uniqueness}, the relative universal projective system is unique up to a unique isomorphism.
\end{rmk}

\begin{prop}[{\cite[Theorem 3.1 and Remark 3.4]{CDPS}}]\label{prop: CDPS universal objects}
There exists a relatively universal projective system $\{(\cW_n, w_n)\}_{n\geq 1}$ in $\Uni_{\dR}(q_0,b_0)$. Moreover, the formation of the relatively universal projective system is compatible with arbitrary base field extension $E/F$.
\end{prop}

Put $B_n=b_0^*\cW_n$. It is equipped with a log connection
\[
  B_n \ra B_n \otimes \Omega_{\ul T/F}^1, \ \ v \mapsto D_n(v)\otimes d \log t
\]
induced from that on $\cW_n$. The nilpotent operator $D_n$ is called the \emph{residue derivation} on $B_n$. 
Taking $\cE= \cW_n$ in \eqref{eq: relative evaluation 2}
induces an isomorphism 
\[
    \ev_{w_n}\colon q_{0,\dR,*}\sEnd(\cW_n)\xra{\sim}B_n.
\]
Then the composition on $q_{0,\dR,*}\sEnd(\cW_n)$ induces an associative $F$-algebra structure on $B_n$. Moreover, for any $r,s \geq 1$, we have
\[
 \uni(\cW_{r} \otimes \cW_s) \leq \uni(\cW_{r})+ \uni(\cW_s)- 1 \leq r+s-1,
\]
and universality induces morphisms
\[
 c_{r,s} \colon (\cW_{r+s-1},w_{r+s-1}) \ra (\cW_{r} \otimes \cW_s, w_r\otimes w_s).
\]
Pulling back along $b_0$, we obtain coproduct maps $c_{r,s} \colon B_{r+s-1} \ra B_{r} \otimes B_s$. Put
\[
 \cH_{\dR} \coloneqq \varinjlim_n B_n^\vee = \varinjlim_n (b_0^*\cW_n)^\vee.
\]
The discussion above makes $\cH_{\dR}$ into a commutative Hopf algebra in $\Ind (\MIC^{\nil}(\ul T/F) )$ with a residue derivation $-D^\vee \coloneqq -\varinjlim_n D_n^\vee$.

\begin{defn}
    The \emph{relative unipotent log de Rham fundamental group} of $(q_0, b_0)$ is defined as 
    \[
    \cG = \pi_1^{\dR}(\ul Y/\ul T, b_0) \coloneqq \sSpec_{\MIC^{\nil}(\ul T/F)}\cH_\dR,
    \]
    which is a pro-unipotent affine group scheme in $\MIC^{\nil}(\ul T/F)$. The \emph{monodromy automorphism} $N$ on $\cG$ is characterized by $N^*=\exp(-D^\vee)$.\footnote{This definition of monodromy agrees with the convention in \cite[\S 6.1]{CDPS}. 
    } Let $G$ denote the underlying group scheme of $\cG$.
\end{defn}

\begin{defn}
    Let $\Rep_{\MIC^\nil(\ul T/ F)}(\cG)$ be the representation category of $\cG$ in $\MIC^\nil(\ul T/ F)$. Equivalently, $\Rep_{\MIC^\nil(\ul T/ F)}(\cG)$ is the category consisting of triples $(V,N_V, \rho_V)$, where $V$ is a finite-dimensional $F$-vector space, $N_V\in\End_F(V)$ is a nilpotent endomorphism, and 
    \[
       \rho_V \colon G\longrightarrow\GL(V)
    \]
    a representation of $G$ satisfying
    \[
       \rho_V\circ N = \mathrm{Ad}_{\exp(-N_V)}\circ \rho_V.
    \]
    Morphisms between $(U,N_U, \rho_U)$ and $(V,N_V, \rho_V)$ are morphisms $\phi\colon U \ra V$ between $G$-representations such that $\phi\circ N_V= N_W\circ \phi$.
\end{defn}

\begin{rmk}
\begin{enumerate}
    \item The proof of Proposition~\ref{prop: CDPS universal objects} uses a construction similar to that of Theorem~\ref{thm: existence of relative universal system}. In particular, $\cW_1= (\cO_{\ul Y}, d)$, and $(\cW_n)_{n\geq 1}$ fits into short exact sequences
    \begin{equation}\label{eq: de Rham extensions}
        0 \ra q_0^*(\mR^1q_{0,\dR,*}\cW_n^\vee)^\vee \ra \cW_{n+1} \ra \cW_n \ra 0.
    \end{equation}
    \item The relatively universal objects and relative unipotent de Rham fundamental group can be constructed for more general pointed relative varieties, see \cite[\S1 (A)-(E)]{CDPS}.
    \item We can also define $\pi_1^{\dR}(\ul Y/ \ul T, b_0)$ to be the Tannakian dual of $(\MIC^{\nil}(\ul Y/\ul T), b_0^*)$. The equivalence between the two constructions is stated in \cite[Theorem 3.12]{CDPS}. In particular, by \emph{loc.\ cit.} and Tannakian duality, there is a natural equivalence of categories
    \begin{equation}\label{eq: Tannakian duality}
            b_0^*\colon \Uni_{\dR}(q_0) \xra[]{\sim} \Rep_{\MIC^{\nil}(\ul T/F)}(\cG).
    \end{equation}
    Under this equivalence, $B_n$ corresponds to $\cW_n$. 
    \item Following \cite[\S 2]{CDPS}, we may define the (absolute) unipotent log de Rham fundamental groups $\pi_1^{\dR}(\ul Y/ F, b_0)$ and $\pi_1^{\dR}(\ul T/ F, b_T)$ to be the Tannakian duals of $(\MIC^{\nil}(\ul Y/ F), b_0^*)$  and $(\MIC^{\nil}(\ul T/ F), b_T^*)$, respectively, where $b_T\colon T \ra T$ denotes the identity section. By \cite[Proposition 2.3]{CDPS}, there is a split exact sequence 
    \begin{equation}\label{eq: split exact sequnence}
    \begin{tikzcd}[column sep=0.8cm]
        1 \ar[r] & \cG \ar[r] &\pi^{\dR}_1(\ul Y/F, b_0)\ar[r, "q_{0,*}"] &\pi^{\dR}_1(\ul T/F, b_T) \ar[r] \ar[l, bend right =30pt, "b_{0,*}"'] & 1.
    \end{tikzcd}
    \end{equation}
    Therefore, $\pi^{\dR}_1(\ul T/F, b_T) \cong \bG_{a,F}$ acts on $\cG$ by conjugation. The monodromy $N$ defined above identifies with the conjugate action of $1\in \bG_{a,F}$ (cf. \cite[\S 6.1]{CDPS}). 
\end{enumerate}
\end{rmk}

For any field extension $E/F$ and a relatively universal system $\{(\cW_n, w_n)\}_{n\geq 1}$ in $\Uni_{\dR}(q_0,b_0)$, its base change $\{(\cW_{n,Y_E}, w_{n,Y_E})\}_{n\geq 1}$ is a relative universal projective system in $\Uni_{\dR}(q_{0,E},b_{0,E})$ (cf. Proposition~\ref{prop: CDPS universal objects}). We put $B_{n,Y_E}= b_{0,E}^*\cW_{n,Y_E}$ and let $\cG_{Y_E}$ be the corresponding relative unipotent log de Rham fundamental group of $(\ul {Y_E}/\ul{T_E},b_{0,E})$.

\begin{lem}[{\cite[Remark 3.3]{CDPS}}]\label{lem: field extension}
For any field extension $E/F$, scalar extension induces an isomorphism
 \[
  \cG_{Y_{E}} \cong  \cG_{Y} \otimes_F E.
 \]
\end{lem}

\vspace{0.05in}
\subsection{Formal deformations}
\label{subsec: formal deformation}
\noindent
\vspace{0.05in}

\noindent 
We have the following criterion for smoothness of log curves in terms of certain commutator quotient.

\begin{thm}\label{thm: HH_0 criterion for smoothness}
    Let $(q_0\colon \ul Y \ra \ul T,b_0)$ be a geometrically connected minimal semistable projective pointed log curve over $F$. Then $Y$ is smooth if and only if $D_5=0$ on $\HH_0(B_5)$.
\end{thm}

Theorem~\ref{thm: HH_0 criterion for smoothness} will be proved in \S \ref{subsec: proof of thm 4.9}. Before that, let us show how it implies Theorem~\ref{thm: HH_0 good reduction}. In the rest of \S \ref{subsec: formal deformation}, let $Z$ be a geometrically connected smooth projective curve over $K$ of genus $g\geq2$ with structure map $q\colon Z\rightarrow \Spec\, K$ and a section $b$. 

Suppose that $Z$ has semistable reduction, with minimal regular model $\cZ$ over $\cO_K$ (cf. \cite[Tag 0C2V]{stacks}). The section $b$ lifts to a \emph{good section} $b_{\OK}$, i.e, a section lying in the relative smooth locus of $\cZ$ over $\cO_K$.\footnote{
By the valuation criterion of properness, $b$ uniquely extends to a section $b_{\cO_K}$ in $\cZ(\cO_K)$. It must be a good section because the model is (strictly) semistable.}
Notice that $\cZ$ is semistable, and $Z$ has good reduction if and only if $\cZ$ is smooth over $\cO_K$. We call $\cZ$ a \emph{minimal semistable model} of $Z$. In particular, $\cZ$ has no exceptional curve of the first kind. Moreover, we equip both $\cZ$ and $\mathcal{O}_K$ with the log structure given by the special fiber, so that the structure morphism $q_{\cO_K}\colon \cZ\rightarrow \mathrm{Spec}\,\cO_K$ is log smooth.

Fix a uniformizer $\varpi$ of $K$, and put
\[
 W=W(k),\ \  K_0=W[\frac 1p],\ \
 \widetilde S=(\Spec W[\![z]\!], \bN, 1\mapsto z)^a.
\]
There are two closed embeddings 
\[
  \rho_\varpi \colon \Spec \cO_K \ra \widetilde{S}, \ \ 
  \rho_0 \colon \Spec W \ra \widetilde{S},
\]
with the corresponding morphisms on the structure rings sending $z$ to $\varpi$ and $0$, respectively.

\begin{construction}[{Formal deformation}]
\label{constr: arithmetic log deformation}
By \cite[\S 2]{AIK} and \cite[\S 6.8]{CDPS}, the pair $(\cZ, b_{\cO_K})$ admits a formal deformation over $\widetilde{S}$; namely, there exists a proper log smooth integral family $\widetilde q:\widetilde{\cZ}\ra\widetilde S$ together with a good section $\widetilde b:\widetilde S\ra\widetilde{\cZ}$ such that
\[
 (\widetilde{\cZ},\widetilde b)
 \times_{\widetilde S,\rho_\varpi}\Spec\OK
 \cong(\cZ,b_{\OK}).
\]
Note that the log structure on $\cZ$ is the one given by the fiber along $\rho_0$. 

Let $R\coloneqq K_0[\![x]\!]$ and $S_R\coloneq  (\Spec R,\bN, 1\mapsto px )^a$.
The composition
\[
 W[\![z]\!]\xra{z\,\mapsto\, px}W[\![x]\!]
 \hookrightarrow R=K_0[\![x]\!]
\]
induces a morphism of log schemes $S_R\rightarrow \widetilde{S}$. Consider the base change
\[
 (\cZ_R,b_R)\coloneq (\widetilde{\cZ},\widetilde b)
 \times_{\widetilde{S}}S_R,
\]
and write $q_R:\cZ_R\ra S_R$ for the structure morphism.

Let $T=\Spec K_0$ and $\ul T\coloneqq  (\Spec K_0,\bN, 1\mapsto 0 )^a$. Consider base change
\[\ul Y\coloneqq \cZ_R\times_{S_R}\ul T\]
and let $b_0\coloneq b_R|_T$. Let $q_0\colon \ul Y \ra \ul T$ denote the structure map. Then $(q_0\colon \ul Y \ra \ul T, b_0)$ is a pointed log curve over $K_0$.
We obtain Cartesian squares of pointed log schemes
\[
\begin{tikzcd}[column sep=1.6cm, row sep=1.2cm]
 \cZ\ar[r]\ar[d, "q_{\cO_K}"]&\widetilde{\cZ}\ar[d,"\widetilde q"]
 & \cZ_R\ar[l]\ar[d, "q_R"']& \ul Y\ar[l]\ar[d, "q_0"']\\
 \Spec \cO_K \ar[u, bend left= 40 pt, "b_{\cO_K}"]\ar[r,"\rho_\varpi"]& \widetilde S \ar[u, bend left= 40 pt, "\widetilde b"{pos=0.5, yshift=-5pt}]
 & S_R \ar[l] \ar[u, bend right= 40 pt, "b_R"'] & \ul T \ar[l] \ar[u, bend right= 40 pt, "b_0"'] .
\end{tikzcd}
\]
\end{construction}

\begin{lem}[{\cite[\S 6.8]{CDPS}}]\label{lem: deformation preserves special graph}
The curve $\ul Y/\ul T$ is a minimal semistable log curve in the sense of Definition~\ref{def: minimal semistable}. Its geometric dual graph\footnote{For the reader's convenience, let us recall the notion of \emph{geometric dual graph} of a nodal curve. Let $\overline F$ be an algebraic closure of $F$. For a nodal curve $Y$ over $F$, the geometric dual graph $\Gamma(Y)$ of $Y$ has one vertex for each irreducible component of $Y_{\overline F}$ and one edge for each of its nodes. The endpoints of an edge correspond to the components containing the two branches at the associated node. If both branches belong to the same component, the edge is a self-loop.} agrees with that of $\cZ_k$. In particular, $Z$ has good reduction if and only if the underlying curve $Y$ of $\ul Y$ is smooth.
\end{lem}

\begin{lem}\label{lem: periods commute with HH_0}
Let $A$ be a finite-dimensional $\Gal_K$-equivariant $\Qp$-algebra whose underlying $\Gal_K$-representation is semistable. Then there is a $K_0$-algebra structure on $D_{\st,K}(A)$ induced from that on $A$. Moreover, $\HH_0(A)$ is semistable and there is a natural isomorphism
\begin{equation}\label{eq: periods commute with HH_0}
 D_{\st,K}(\HH_0(A))
 \cong\HH_0(D_{\st,K}(A)),
\end{equation}
compatible with monodromy.
\end{lem}

\begin{proof}
First, since $[A,A]$ is a subrepresentation of $A$, we know that $\HH_0(A)=A/[A, A]$ is also a semistable representation. Secondly, since $A$ is semistable, there is an isomorphism between filtered $(\vp,N)$-modules 
\[
 D_{\st,K}(A\otimes_{\Qp}A) \cong D_{\st,K}(A)\otimes_{K_0}D_{\st,K}(A).
\]
The multiplication on $A$ then induces a $K_0$-algebra structure on $D_{\st,K}(A)$.
Applying $D_{\st,K}$ to the exact sequence
\[
 A\otimes_{\Qp} A \xra{c_A} A \ra \HH_0(A) \ra 0, \ \ c_A(a_1\otimes a_2)= a_1a_2-a_2a_1 
\]
yields an exact sequence
\[
 D_{\st,K}(A)\otimes_{K_0} D_{\st,K}(A) \ra D_{\st,K}(A) \ra D_{\st,K}(\HH_0(A)) \ra 0
\]
in the category of filtered $(\vp,N)$-modules. Note that the first arrow is the commutator map with respect to the induced $K_0$-algebra structure; taking cokernel induces \eqref{eq: periods commute with HH_0}. Since arrows commute with monodromy, so does \eqref{eq: periods commute with HH_0}. 
\end{proof}

Let $\{(\bE_{n,Z},u_{n,Z})\}_{n\geq 1}$ be a relatively universal projective system in $\Uni_{\Zp}(q,b)$. For the pointed log curve $(q_0\colon \ul Y\rightarrow \ul T, b_0)$, let $\{(\cW_{n,Y},w_{n,Y})\}_{ n\geq 1}$ be a relatively universal projective system in $\Uni_{\mathrm{dR}}(q_0, b_0)$. Put
\[
 A_{n,Z}=b^*\bE_{n,Z}, \ \ B_{n,Y}=b_0^*\cW_{n,Y}.
\]

\begin{prop}\label{prop: semistable comparison}
    There is an isomorphism
    \begin{equation}\label{eq: semistable comparison}
        D_{\st, K}(A_{n,Z}[\frac 1p]) \xra{\sim} B_{n,Y}, 
    \end{equation}
    of augmented $K_0$-algebras, compatible with monodromy automorphisms, where the monodromy on $B_{n,Y}$ is the residue derivation $D_n$ of its log connection. In particular, it induces an isomorphism 
    \begin{equation}\label{eq: semistable comparison after HH_0}
          D_{\st, K} (\HH_0 (A_{n,Z})[\frac 1p] ) \xra{\sim} \HH_0(B_{n,Y}).  
    \end{equation}
\end{prop}

\begin{proof}
    The isomorphism \eqref{eq: semistable comparison} is constructed in \cite[\S 6.8]{CDPS}, 
    while the isomorphism \eqref{eq: semistable comparison after HH_0} follows from Lemma~\ref{lem: periods commute with HH_0}.
\end{proof}

\begin{prop}
    Theorem~\ref{thm: HH_0 criterion for smoothness} implies Theorem~\ref{thm: HH_0 good reduction}.
\end{prop}

\begin{proof}
    Assume that Theorem~\ref{thm: HH_0 criterion for smoothness} holds. By Lemma~\ref{lem: deformation preserves special graph}, $Z$ has good reduction over $\mathcal{O}_K$ if and only if $Y$ is smooth. By Theorem~\ref{thm: HH_0 criterion for smoothness}, this is equivalent to saying that $D_5=0$ on $\HH_0(B_{5,Y})$. It is also equivalent to the monodromy being zero on $D_{\st,K}(\HH_0(A_{5,Z})[\frac 1p])$ by Proposition~\ref{prop: semistable comparison}, i.e., $\HH_0(A_{5,Z})[\frac 1p]$ is crystalline. 
\end{proof}

\vspace{0.05in}
\subsection{Proof of Theorem~\ref{thm: HH_0 criterion for smoothness}}
\label{subsec: proof of thm 4.9}
\noindent
\vspace{0.05in}

\noindent 
Lastly, we complete the proof of Theorem~\ref{thm: HH_0 criterion for smoothness}. Lemma~\ref{lem: field extension} allows us to enlarge the base field. Therefore, we may assume (throughout \S \ref{subsec: proof of thm 4.9}) that $F$ is an algebraically closed field of characteristic $0$. Let $(q_0\colon \ul Y \ra \ul T, b_0)$ be a pointed log curve over $F$ and let $\cG= \pi_1^{\dR}(\ul Y/\ul T, b_0)$ be the relative unipotent log de Rham fundamental group of $(q_0,b_0)$. Let $G$ denote the underlying group scheme of $\cG$.
\begin{defn}
We define the \emph{pro-nilpotent de Rham fundamental Lie algebra} of $(\ul Y/ \ul T,b_0)$  to be $\kg= \kg_Y \coloneqq \Lie(G)$.
Let $D_{\kg}= -\log(dN)$ be the \emph{residue derivation} on $\kg$ induced from the monodromy $N$ on $\cG$.\footnote{Recall that our monodromy is defined as $N^*= \exp(-D^\vee)$.
}
\end{defn} 
Let $U(\kg)$ be the universal enveloping algebra of $\kg$. Let $I \subseteq U(\kg)$ be the two-sided ideal generated by $\kg$. Then $U(\kg)$ admits an augmentation $\varepsilon\colon U(\kg)\ra F$ given by $\varepsilon(1)=1$ and $\varepsilon(x)=0$ for all $x\in\kg$, with $\Ker\varepsilon = I$. The derivation $D_\kg$ extends to derivations $D_{U(\kg)}$ on $U(\kg)$ and $D_{U_n}$ on $U_n \coloneqq U(\kg)/I^n$.

\begin{lem} For any $n \geq 1$, 
    \begin{enumerate}
        \item $D_{U_n}$ is nilpotent on $U(\kg)/I^n$;
        \item $U(\kg)/I^n$ is finite-dimensional over $F$.
    \end{enumerate}
\end{lem}
\begin{proof}
    Consider the lower central series
    \[
     C^1\kg= \kg,\ \  C^{j+1}\kg= [\kg, C^j\kg], \ \ \forall j \geq 1,
    \]
    and put $\kg_n= \kg/C^n\kg$, which is nilpotent with a nilpotent derivation $D_{\kg_n}$ induced from that on $\kg$. Let $J_n \subseteq U(\kg_n)$ be the ideal generated by $\kg_n$. Since $C^n\kg \subseteq I^n$, we have a natural isomorphism 
    \[
     U(\kg)/I^n \cong U(\kg_n)/J_n^n.
    \]
    Since $D_{\kg_n}$ is nilpotent, the induced derivation on $U(\kg_n)/J_n^n$ is also nilpotent. This proves (1). 
    For (2), note that $\kg/C^2\kg \cong \mH_\dR^1(\ul Y/\ul T)$ is finite-dimensional over $F$. The surjection $\kg \ra \kg_n$ ensures that $\kg_n/C^2\kg_n$ is also finite-dimensional over $F$. Therefore, since $\kg_n$ is nilpotent, $\kg_n$ is also finitely generated.\footnote{Indeed, let $\kh$ be a finitely generated sub-Lie-algebra of $\kg_n$ such that $\kg_n \subseteq \kh+ C^2\kg_n$, one can use induction to show that $\kg_n \subseteq \kh+ C^j\kg_n$ for all $j \geq 1$.} This forces $U(\kg_n)/J_n^n$ to be finite-dimensional.
\end{proof}

\begin{defn}
  We define $\mathrm{Mod}^{\mathrm{fd}}_{U(\kg)}$ to be the category of  finite-dimensional $F$-vector spaces $V$ together with a nilpotent endomorphism $D_V$ and a left continuous $U(\kg)$-module structure such that 
 \begin{equation}\label{eq: Leibniz rule} 
   D_V(av)= D_{U(\kg)}(a)v+a D_V(v), \ \ \forall a \in U(\kg), \ v\in V.   
 \end{equation}
 The continuity is with respect to the pro-nilpotent topology on  $U(\kg)$. 
 Morphisms are morphisms of $U(\kg)$-modules compatible with their derivations. Similarly, the category $\mathrm{Mod}^{\mathrm{fd}}_{U(\kg)/I^n}$ is defined in a similar way. Moreover, let $\mathrm{Mod}^{\mathrm{fd},\leq n}_{U(\kg)}$ denote the full subcategory of $\mathrm{Mod}^{\mathrm{fd}}_{U(\kg)}$ consisting of objects $V$ with $I^n V =0$. Note that by definition, $\mathrm{Mod}^{\mathrm{fd},\leq n}_{U(\kg)}$ is naturally equivalent to $\mathrm{Mod}^{\mathrm{fd}}_{U(\kg)/I^n}$.
\end{defn}

On the $\cG$-representation side, let $\Rep^{\leq n}_{\MIC^{\nil}(\ul Y/F)} (\cG)$ be the full subcategory of $\Rep_{\MIC^{\nil}(\ul Y/F)} (\cG)$ consisting of object $V$ that admits a filtration
 \[
V= V^1 \supseteq V^2 \supseteq \cdots \supseteq V^n \supseteq V^{n+1}=0
 \]
such that $V^i/V^{i+1}$ is a trivial $\cG$-representation for each $i \geq 1$. Under Tannakian duality \eqref{eq: Tannakian duality}, there is a natural equivalence
 \[
  b_0^*\colon \Uni_{\dR}^{\leq n}(\ul Y/\ul T) \xra[]{\sim} \Rep^{\leq n}_{\MIC^{\nil}(\ul T /F)} (\cG).
 \]

Consider the natural functor
 \begin{equation}\label{eq: G-representation to U(g)-module}
     \mathrm{Diff}_\cG \colon \Rep_{\MIC^{\nil}(\ul T/F)} (\cG) \ra \mathrm{Mod}^{\mathrm{fd}}_{U(\kg)}
 \end{equation}
induce by differentiation.

\begin{lem}\label{lem: categorical equivalence}
The functor \eqref{eq: G-representation to U(g)-module} induces a natural equivalence of categories
 \[
  \Rep^{\leq n}_{\MIC^{\nil}(\ul T/F)} (\cG) \xra{\sim} \mathrm{Mod}^{\mathrm{fd},\leq n}_{U(\kg)} \cong \mathrm{Mod}^{\mathrm{fd}}_{U(\kg)/I^n}.
 \]
\end{lem}
\begin{proof}
First note that $\mathrm{Diff}_\cG$ is fully faithful. Indeed, this holds for unipotent groups by \cite[Chapter 10, \S I]{Milne2017}, and any finite-dimensional representation of a pro-unipotent group factors through a finite dimensional unipotent quotient. 
Secondly, we show that $\mathrm{Diff}_\cG$ induces a fully faithful functor 
\[
 \mathrm{Diff}_\cG^{\leq n} \colon \Rep^{\leq n}_{\MIC^{\nil}(\ul T/F)} (\cG) \ra \mathrm{Mod}^{\mathrm{fd},\leq n}_{U(\kg)}
\]
Take any $V \in \Rep^{\leq n}_{\MIC^{\nil}(\ul Y/\ul T)} (\cG)$ with filtration
 \[
V= V^1 \supseteq V^2 \supseteq \cdots \supseteq V^n \supseteq V^{n+1}=0
 \]
such that $V^i/V^{i+1}$ is a trivial representation of $\cG$ for each $i$, then $V^i/V^{i+1}$ is also a trivial $U(\kg)$-module. Thus, we must have $IV^{i} \subseteq V^{i+1}$. This implies $I^n V \subseteq V^{n+1}=0$. The full faithfulness of $\mathrm{Diff}_\cG$ implies the full faithfulness of $\mathrm{Diff}_\cG^{\leq n}$. 

We are left to check that $\mathrm{Diff}_\cG^{\leq n}$ is essentially surjective. Let $\exp_\kg \colon \kg \ra G$ be the exponential. Note that $\exp_\kg$ is an isomorphism on the underlying schemes with inverse $\log_G\colon G \ra \kg$.  For any $(V, D_V)\in \mathrm{Mod}^{\mathrm{fd},\leq n}_{U(\kg)}$, define $\rho_V \colon G \ra \GL(V)$ by
\[
\rho_V(\exp_\kg x)\coloneqq \exp(L_x)= \sum_{r=0}^{n-1}\frac{(L_{x})^r}{r!},
\]
for all $x\in \kg$, where $L_x$ denotes the action of $x$ on $V$. 
We claim that $\rho_V$ is compatible with derivation. Denote $\widetilde D_{\kg} \coloneqq \exp(-D_\kg)$ and $\widetilde D_{V} \coloneqq \exp(-D_V)$. For any $x \in \kg$, the Leibniz rule \eqref{eq: Leibniz rule} induces $L_{D_\kg(x)}= [{D_V},L_{x}]$. Recall that $dN= \exp(-D_\kg)= \widetilde D_\kg$, so we have
\begin{align*}
    \rho_V  (N(\exp_\kg x) )= & \rho_V (\exp_\kg (d N(x)) )= \rho_V (\exp_\kg (\widetilde D_\kg(x)) ) \\
    =& \exp ({L_{\widetilde D_\kg(x)}} ) 
    =  \exp (\mathrm{Ad}_{\widetilde D_V}L_x )= \mathrm{Ad}_{\widetilde D_V} ((\rho_V)(\exp_\kg x) ).
\end{align*}
So $(\rho_V, D_V)$ is an object in $\Rep_{\MIC^{\nil}(\ul T/ F)} (\cG)$. By construction, $\mathrm{Diff}_\cG(\rho_V)= V$. Now we check the unipotency index of $V$. For any $i\geq 1$, take $V^i= I^{i-1}V$, then $V^{n+1}=0$. Note that 
\[
 \rho_V(x)(v)= \sum_{r=0}^{n-1}\frac{(\log_G x)^r v }{r!} \in V^i,
\]
for all $v\in V^i$ and $x\in G$.
Moreover, using $D_{U(\kg)}(I) \subseteq I$ and the Leibniz rule, we have $D_V(V^i) \subseteq V^{i}$. Thus $V^i$ is a $\cG$-subrepresentation of $V$. Finally, using the fact $I(V^i/V^{i+1})=0$ and the full faithfulness of $\mathrm{Diff}_\cG$, we conclude that $(\rho_V,D_V)$ is an object in $\Rep^{\leq n}_{\MIC^{\nil}(\ul T/F)} (\cG)$. This finishes the proof.
\end{proof}

Let $\{(\cW_n, w_n)\}_{n\geq 1}$ be a relatively universal projective system in $\Uni_{\dR}(q_0, b_0)$. Let $B_n= b_0^*\cW_n$ with residue derivation $D_n$. The Tannakian duality \eqref{eq: Tannakian duality} sends $\cW_n$ to $B_n$. In particular, $B_n$ is a universal object in $\Rep^{\leq n}_{\MIC^{\nil}(\ul T /F)} (\cG)$. Combining this with Lemma~\ref{lem: categorical equivalence}, we obtain:

\begin{prop}\label{prop: Lie algbebra vs B_n}
There are canonical isomorphisms 
\[
 B_n^{\mathrm{op}}\cong U(\kg)/I^n
\]
of augmented $F$-algebras in $\MIC^{\nil}(\ul T/ F)$, compatible with transition maps.
Under these isomorphisms, the residue derivation $D_n$ on $B_n$ coincides with the residue derivation induced by $D_{\kg}$.
\end{prop}

\begin{proof}
By Lemma~\ref{lem: categorical equivalence} and the fact that $(B_n,1_n=w_n(1))$ is universal in $\Rep^{\leq n}_{\MIC^{\nil}(\ul T /F)} (\cG)$, there is a unique morphism 
\[
\beta \colon B_n \ra U(\kg)/I^n,
\]
sending unit to unit, and compatible with derivations. It must be an isomorphism because $U(\kg)/I^n$ shares the same universality in $\mathrm{Mod}^{\mathrm{fd}}_{U(\kg)/I^n}$. The compatibility of augmentation and transition maps also follows from universality. Therefore, we are left to check that $\beta$ is compatible with multiplication. If we view $B_n$ as a $U(\kg)/I^n$-module and consider
\[
  \gamma \colon U(\kg)/I^n \ra B_n, \ \ a \ra a\cdot 1_n,
\]
then $\beta$ and $\gamma$ are inverse to each other. For any $x\in B_n$, let $T_x$ be the unique endomorphism of $B_n$ sending $1_n=w_n(1)$ to $x$. It must be a $U(\kg)/I^n$-morphism by uniqueness. The multiplication on $B_n$ is therefore given by $xy=  T_x(y)$. For any $a, c\in U(\kg)/I^n$, we have 
\[
 \gamma(a)\gamma(c)= T_{\gamma(a)}( \gamma (c))= T_{\gamma(a)}(c\cdot 1_n)= c\cdot T_{\gamma(a)}( 1_n)= ca\cdot 1_n= \gamma(ca).
\]
as desired.
\end{proof}

We provide an explicit description of the pro-nilpotent de Rham fundamental Lie algebra in the case where the dual graph of the log curve is a chain. For a finite-dimensional $F$-vector space $V$, let $\kl(V)$ be the free pro-nilopotent Lie algebra generated by $V$, and let $L(V) \coloneqq \bigoplus_{d\geq 0}V^{\otimes d}$.

\begin{prop}[{\cite[Remark 6.28]{CDPS}}]\label{prop: chain Lie presentation}
Let $\ul W/ \ul T$ be a minimal semistable pointed log curve over $F$. Assume that $W$ has components $W_0,W_1,\ldots,W_s$ for some $s\geq 1$ such that its dual graph is a chain 
\[ 
  W_0-W_1-\cdots -W_s,
\]
and such that $W_i$ is isomorphic to $\bP^1_F$ for every $0<i<s$. Take $c_0\in W_0\cap W_1$ as a section of $W$. Let $g_L$ and $g_R$ be the genera of $W_0$ and $W_s$, respectively, both of which are positive by the minimality of $W$. Let $V_L= \mathrm{span}_F\{a_1,b_1,\ldots,a_{g_L},b_{g_L}\}$ and $V_R= \mathrm{span}_F\{c_1,d_1,\ldots,c_{g_R},d_{g_R}\}$ be two vector spaces of dimension $2g_L$ and $2g_R$, respectively. Then there is an isomorphism 
\[
 \kg_W\cong
 \kl(V_L\oplus V_R)\big/\langle e_L-e_R\rangle.
\]
of pro-nilpotent Lie algebras, where
\[
  e_L =\sum_{i=1}^{g_L}[a_i,b_i], \ \  e_R=-\sum_{j=1}^{g_R}[c_j,d_j].
\]
Moreover, let $\{(\cW_{n,W},w_{n,W})\}_{n\geq 1}$ be a relatively universal projective system in $\Uni_{\dR}(\ul W/\ul T, c_0)$ and let $B_{n,W}=c_0^*\cW_{n,W}$. Then 
\begin{equation}\label{eq: chain associative presentation}
 B_{n,W}\cong
 L(V_L\oplus V_R)
 \big/ (J^n+(e_L-e_R) ),
\end{equation}
where $J=\bigoplus_{d\geq1}(V_L\oplus V_R)^{\otimes d}$ and $(e_L-e_R)$ is the two-sided ideal generated by $e_L-e_R$.
\end{prop}

\begin{proof}
This follows from Proposition~\ref{prop: Lie algbebra vs B_n} and \cite[Remark 6.28]{CDPS}.
\end{proof}

\begin{prop}\label{prop: chain monodromy derivation}
Let $W$ be as above. Then there is a $\lambda\in F^\times$ such that the derivation $D_\kg$ on $\kg_W$ (and hence the residue derivation $D_n$ on every $B_{n,W}$) satisfies
\begin{equation}\label{eq: chain Lie algebra derivation}
D_\kg(v)= 
\begin{cases}
    0 \ \ &v\in V_L, \\
   \lambda[e,v]\ \ &v\in V_R,
\end{cases}
\end{equation}
where $e$ denotes the common image of $e_L$ and $e_R$ in $\kg_{W}$.
\end{prop}

\begin{proof}
    Let $W_{\geq 1}\coloneqq W_1\cup\cdots \cup W_s$. Define 
    \[
     \cG^L \coloneqq \mathrm{image} (\pi_1^\dR(\ul {W_0}/\ul T,c_0) \ra \cG_{W} ), \ \ \cG^R \coloneqq \mathrm{image} (\pi_1^\dR(\ul {W_{\geq 1}}/\ul T,c_0) \ra \cG_{W} )
    \]
    Let $\kg^L$ and $\kg^R$ be the Lie algebras of $\cG^L$ and $\cG^R$, respectively.
    Let $\eta$ be the image of a nonzero element of $\pi_1^\dR(\ul T/ F, T) \cong \bG_{a,F}$ in $\cG_{W}$ under $c_{0,*}$. By \cite[Remark 6.28]{CDPS}, we may write $\eta= \exp_{\kg}(\mu e)$ for some $\mu \in F^\times$. Then by \cite[Proposition 6.27]{CDPS}, the monodromy $N$ satisfies
    \[
    \begin{cases}
      N(\alpha)= \alpha, \ \ &\alpha\in \cG^L, \\ 
      N(\beta)= \eta^{-s} \beta \eta^s, \ \ & \beta\in \cG^R.  
    \end{cases}
    \]
    Recall that $D_{\kg}= -\log dN$.
    Thus, by differentiating the group automorphism $N$, the induced action on $\kg_W$ is given by 
    \begin{equation}\label{eq: monodromy formula}
        d N|_{\kg^L}= \id_{\kg^L}, \ \ d N|_{\kg^R}= \exp(-s\mu \cdot\mathrm{ad}_e).
    \end{equation}
    Taking the logarithm of \eqref{eq: monodromy formula} yields \eqref{eq: chain Lie algebra derivation} with $\lambda=  s\mu$. 
\end{proof}

\begin{lem}\label{lem: chain HH_0 trace}
For $W$ as above, the derivation $D_5$ is nonzero on $\HH_0(B_{5,W})$.
\end{lem}

\begin{proof}
Keep the notation as in Proposition~\ref{prop: chain Lie presentation} and put $L\coloneqq  L(V_L\oplus V_R)$. Let $R\coloneqq F[\tau]/(\tau^5)$ and denote
\[
 P\coloneqq \begin{pmatrix}0&1\\0&0\end{pmatrix},\ \ 
 Q\coloneqq \begin{pmatrix}0&0\\1&0\end{pmatrix},\ \ 
 H\coloneqq [P,Q]=\begin{pmatrix}1&0\\0&-1\end{pmatrix}.
\]
First, consider the homomorphism $\widetilde\rho\colon L\ra M_2(R)$ of unital $F$-algebras defined by
\[
 \widetilde\rho(a_1)=\widetilde\rho(d_1)= \tau P,\ \  \widetilde\rho(b_1)= \widetilde\rho(c_1)= \tau Q,
\]
while sending all the other generators to zero.
Notice that 
\[
[\tau P,\tau Q]=-[\tau Q,\tau P]=\tau^2H,
\]
this shows that $\widetilde\rho(e_L-e_R)=0$.
Moreover, every generator has image belonging to $\tau M_2(R)$, so
\[
 \widetilde\rho(J^5)\subseteq\tau^5M_2(R)=0.
\]
Therefore, $\widetilde\rho$ factors through a well-defined $F$-algebra homomorphism
\[
 \rho:B_{5,W}\ra M_2(R).
\]

It follows from Proposition~\ref{prop: chain monodromy derivation} that $D_5(a_1c_1)=\lambda a_1[e,c_1]$ in $\HH_0(B_{5,W})$.
Consider
\[
 \ell \coloneqq \tr\circ \rho \colon B_{5,W} \ra R,
\]
where the trace is taken in $M_2(R)$. Since the ring $R$ is commutative, $\ell$ factors through $\HH_0(B_{5,W})$. Since $[H,Q]=-2Q$ and $\tr(PQ)=1$, we have
\begin{equation}\label{eq: nonzero cyclic trace}
 \ell(D_5(a_1c_1))= \tr \rho(\lambda a_1[e,c_1])= 
 \lambda\tau^4\tr (P[H,Q] )
 =-2\lambda\tau^4\neq0.
\end{equation}
Therefore, $D_5(a_1c_1)\neq 0$ in $\HH_0(B_{5,W})$, as desired. 
\end{proof}

\begin{prop}\label{prop: contraction to chain}
Let $(\ul Y/\ul T, b_0)$ be a geometrically connected minimal semistable projective pointed log curve over $F$ . Assume that $Y$ is not smooth.
\begin{enumerate}
    \item If the dual graph of $Y$ contains a loop, then $D_2 \neq 0$ on $B_{2,Y}$.
\item If the dual graph of $Y$ is a tree, then there is a minimal semistable log curve $(\ul W/ \ul T, c_0)$ as in Proposition~\ref{prop: chain Lie presentation}
such that there is a surjection $\kg_Y \twoheadrightarrow \kg_W$ between their pro-nilpotent de Rham fundamental Lie algebras, compatible with derivations up to additions by inner derivations\footnote{
An \emph{inner derivation} is a derivation of the form $x \mapsto [a,x]$.
}.
\end{enumerate}
\end{prop}

\begin{proof}
For (1), the proof of \cite[Theorem 6.24]{CDPS} indicates that the monodromy acts nontrivially on $\Lie(G_Y^{\ab}) \cong \mH_\dR^1(\ul Y/ \ul T)$. The desired nonvanishing then follows from $B_2(Y)\cong F\oplus (\mH_\dR^1(\ul Y/ \ul T))^\vee$, by \eqref{eq: de Rham extensions}. 

Part (2) is already proved in \cite[\S 6.5]{CDPS}. We briefly recall the construction here. Let 
\[
 W_0 - W_1 - \cdots - W_r
\]
be a maximal chain in the dual graph of $Y$ with $r \geq 1$. Since $Y$ is minimal semistable, neither $W_0$ or $W_r$ is isomorphic to $\bP_F^1$. Let $s\geq 1$ be the smallest index such that $W_s$ is not isomorphic to $\bP_F^1$ and take
\[
 W= W_0 \cup W_1 \cup \cdots \cup W_s.
\]
Then the curve $W$, together with the section $c_0\in W_0\cap W_1$, satisfies the desired property. Inner derivations appear when changing the base points, see \cite[Proposition 6.15]{CDPS} for details.
\end{proof}

\begin{proof}[Proof of Theorem~\ref{thm: HH_0 criterion for smoothness}]
    If $Y$ is smooth, $N$ is trivial on $\cG= \pi_1^\dR(\ul Y/\ul T,b_0)$ by \cite[Proposition 6.8(1)]{CDPS}. Therefore, $D_5=0$ on $B_{5,Y}$. Conversely, assume that $Y$ is not smooth. If the dual graph of $Y$ contains a loop, then by Proposition~\ref{prop: contraction to chain} (1), $D_2 \neq 0$ on $B_{2,Y}= \HH_0(B_{2,Y})$. The transition map induces a surjection $B_{5,Y} \twoheadrightarrow B_{2,Y}$, compatible with derivations. Hence $D_5$ is non-zero on $\HH_0(B_{5,Y})$. If the dual graph of $Y$ is a tree, take $W$ as in Proposition~\ref{prop: contraction to chain} (2). 
    Note that inner derivations act trivially on $\HH_0(B_{n,W})$ for each $n\geq 1$. Therefore, The surjection $\kg_Y \twoheadrightarrow \kg_W$ induces a surjection 
    \[
      \HH_0(B_{5,Y}) \twoheadrightarrow \HH_0(B_{5,W}),
    \] 
    compatible with derivations by Proposition~\ref{prop: Lie algbebra vs B_n}. Consequently, we only need to check $D_5\neq 0$ on $\HH_0(B_{5,W})$. This is already handled in Lemma~\ref{lem: chain HH_0 trace}.
\end{proof}

\vspace{0.2in}
\section{Proof of the rigidity theorems}
\label{sec: rigidity theorems}
In this section, we complete the proofs of the main theorems.

\vspace{0.05in}
\subsection{A diagonal construction for a family of curves}
\noindent
\vspace{0.05in}

\noindent 
First, observe that Construction~\ref{constr: diagonal universal system} remains valid in the relative setting. This construction allows us to treat the situation where the family of curves does not admit a section, as mentioned in \S \ref{subsec: intro_proof_main_thm_general}.

\begin{construction}\label{constr: relative diagonal universal system} 
Let $S$ be a smooth rigid analytic variety over $K$ and let $q:X\ra S$ be a smooth proper family of geometrically connected curves of genus $g\geq 2$. 
\begin{enumerate}
\item Consider morphisms
\[
 \pr_1:X\times_S X\ra X,\ \
 \Delta:X\ra X\times_S X,
\]
where $\pr_1$ stands for the projection to the first coordinate, while $\Delta$ is the diagonal. The pair $(\pr_1, \Delta)$ defines a family of point curves as in the setting of \S \ref{sec: universal unipotent local systems}. Applying Theorem~\ref{thm: existence of relative universal system} to  $(\pr_1, \Delta)$,
we obtain a relatively universal projective system $\{(\bE_{n,X}^\Delta, u_{n,X}^\Delta)\}_{n\geq 1}$ in $\Uni_{\Zp}(\pr_1, \Delta)$. Put
\[
 \bA_{n,X}^\Delta\coloneqq \Delta^*\bE_{n,X}^\Delta, \ \ \bL_n\coloneqq  q_* ((\bA_{n,X}^{\Delta})^\vee ).
\]
\item For every classical point $s\in S$ with residue field $\kappa(s)$, let $q_s\colon X_s \ra s$ denote the fiber over $s$ and let $\iota_s \colon X_s \ra X$ denote the closed embedding. Consider morphisms
\[
 \pr_{1,s}:X_s\times_{\Spec\,\kappa(s)} X_s\ra X_s,\ \
 \Delta_s:X_s\ra X_s\times_{\Spec\,\kappa(s)} X_s,
\]
as in (1). Applying Construction~\ref{constr: diagonal universal system} to the pair $(\pr_{1,s}, \Delta_s)$, we obtain a relatively universal projective system $\{(\bE_{n,s}^\Delta, u_{n,s}^\Delta)\}_{n\geq 1}$ in $\Uni_{\Zp}(\pr_{1,s}, \Delta_s)$. Put
\[
 \bA_{n,s}^\Delta\coloneqq \Delta_s^*\bE_{n,s}^\Delta, \ \ C_n(X_s)\coloneqq  \mH^0_\et (X_{\bar s},(\bA_{n,s}^\Delta)^\vee )[\frac 1p],
\]
where ${\bar s}$ is a geometric point above $s$.

\end{enumerate}
\end{construction}

\begin{prop}\label{prop: relative diagonal detector}
Let $X$, $s$ and ${\bar s}$ be as above. For any $n\geq 1$, there is a unique isomorphism 
\begin{equation}\label{eq: base change of diagonal}
\iota_s^*\bA_{n,X}^\Delta \xra{\sim} \bA_{n,s}^\Delta
\end{equation}
of \'etale $\Zp$-local systems on $X_s$, which induces an isomorphism
\begin{equation}\label{eq: base change of diagonal cohomology}
  (\bL_n)_{\bar s}[\frac{1}{p}]\cong C_n(X_s)
\end{equation}
between $\Gal_{\kappa(s)}$-representations.
\end{prop}

\begin{proof}
    Applying Proposition~\ref{prop: base change universal system} to the diagram 
    \[
    \begin{tikzcd}[column sep=1.6cm, row sep=1.2cm]
     X_s\times_{\kappa(s)} X_s \ar[r,"\iota_s \times \iota_s"]\ar[d,"\pr_{1,s}"]& X\times_S X\ar[d,"\pr_1"']\\
     X_s \ar[u, bend left=40pt, "\Delta_s"]
    \ar[r,"\iota_s"]& X\ar[u, bend right=40pt, "\Delta"']
    \end{tikzcd}
    \]
    and pulling back along the diagonals, we obtain \eqref{eq: base change of diagonal}. The isomorphism \eqref{eq: base change of diagonal cohomology} follows from \eqref{eq: base change of diagonal} by taking proper base change.
\end{proof}

\vspace{0.05in}
\subsection{Proof of the main theorems}
\label{subsec: proofs of the main theorems}
\noindent
\vspace{0.05in}

\noindent
Finally, we finish the proofs of Theorem~\ref{thm: introduction rigidity} and Theorem~\ref{thm: introduction rigidity(semistable reduction)}.

\begin{proof}[Proof of Theorem~\ref{thm: introduction rigidity}]
The case for $g=1$ is already handled in \cite[Corollary 1.2]{DY}. In the following, we assume $g \geq 2$.
Consider the local system 
\[
\bL_n=  q_* ((\bA_{n,X}^{\Delta})^\vee ).
\]
By Corollary~\ref{cor: good reduction implies crystalline}, $(\bL_5)_{{\bar s}}[\frac 1p] \cong C_5(X_s)$ is semistable for every $s$. By \cite[Theorem 7.2]{Guoyang}, $\bL_5$ must be a semistable $\Zp$-local system. Moreover, by Corollary~\ref{cor: good reduction implies crystalline} again, $(\bL_5)_{{\bar s}_0}[\frac 1p] \cong C_5(X_{s_0})$ is crystalline. Therefore, \cite[Theorem 1.1]{DY} implies that $(\bL_5)_{\bar s}$ is crystalline for every classical point $s$. Finally, applying Theorem~\ref{thm: diagonal good reduction} completes the proof.
\end{proof}

\begin{proof}[Proof of Theorem~\ref{thm: introduction rigidity(semistable reduction)}]
The case for $g=1$ can be established by applying \cite[Theorem 3.11]{DY} to $\mR^1q_*\ul\Zp$. In the following, we assume $g \geq 2$.
Consider the local system 
\[
\bL_n=  q_* ((\bA_{n,X}^{\Delta})^\vee ),
\]
again. By Corollary~\ref{cor: good reduction implies crystalline} and \cite[Theorem 7.2]{Guoyang}, $\bL_5$ is a semistable $\Zp$-local system. Moreover, by Corollary~\ref{cor: good reduction implies crystalline} again, $(\bL_5)_{{\bar s}_T}[\frac 1p] \cong C_5(X_{s_T})$ is crystalline for every irreducible component $T$ of $\kS_k$.

Then according to \cite[Theorem 3.11]{DY}, $(\bL_5)_{\bar s}[\frac 1p] \cong C_5(X_s)$ is crystalline for every closed point $s$. The desired statement follows from Theorem~\ref{thm: diagonal good reduction} again.
\end{proof}

\vspace{0.05in}
\subsection{Failure of rigidity over a closed annulus}
\label{subsec: failure of rigidity}
\noindent
\vspace{0.05in}

\noindent The following example shows that the good reduction hypothesis on the base in
Theorem~\ref{thm: introduction rigidity} is necessary. 

\begin{example}\label{example: failure of rigidity}
Let
$K=\mathbb{Q}_3$ and let
\[
 S=\operatorname{Sp}(K\langle t,u\rangle/(tu-3))
\]
be the closed annulus $|3|\leq |t|\leq 1$.
The annulus $S$ does not admit a smooth formal model over $\mathcal{O}_K$.
Let $q\colon E\to S$ be the family of projective plane cubics
\[
 Y^2Z=X^3+X^2Z+tZ^3.
\]
Its Weierstrass invariants are $c_4=16$ and $\Delta=-16t(4+27t)$.
Since $t$ is invertible on $S$ and $|27t|<1$, the discriminant is
nowhere zero. Thus $q$ is a smooth proper family of elliptic curves.

For any classical point $s\in S$, consider $t(s) \in \kappa(s)$.
If $|t(s)|=1$, then $\Delta(t(s))$
is a unit in $\mathcal{O}_{\kappa(s)}$, so $E_s$ has good reduction. If $|t(s)|<1$, the special fiber of the Weierstrass model with respect to the chart $Z=1$ is given by 
\[
 y^2=x^2(x+1).
\]
The origin is a node with distinct tangent lines $y=x$ and $y=-x$,
both defined over the residue field of  $\kappa(s)$. Hence $E_s$ has split
multiplicative reduction over $\kappa(s)$; in particular, $E_s$ has semistable reduction.
\end{example}

\vspace{0.2in}

\bibliographystyle{amsalpha}
\raggedbottom
\bibliography{refs}

\vspace{15mm}

\begin{tabular}{l}
    Department of Mathematical Sciences,\\
    Tsinghua University\\
    Beijing, China\\
    \textit{E-mail address: }\texttt{zhefanduan@126.com}
\end{tabular}

\end{document}